\documentclass[12pt]{article}
 \usepackage{latexsym}
 \usepackage{amssymb}
 \usepackage{graphicx}

\newtheorem{Theorem}{Theorem}[part]
\newtheorem{Conjecture}{Conjecture}[part]
\newtheorem{Definition}{Definition}[part]
\newtheorem{Proposition}{Proposition}[part]
\newtheorem{Assumption}{Assumption}[part]
\newtheorem{Lemma}{Lemma}[part]
\newtheorem{Corollary}{Corollary}[part]
\newtheorem{Remark}{Remark}[part]

\def \ep{\hbox{ }\hfill$\Box$}

\def \ra{\rightarrow}
\def\reff#1{{\rm(\ref{#1})}}

\begin{document}
\title{E-Determinants of Tensors}

\author{
Shenglong Hu \thanks{Email: shenglong@tju.edu.cn. Department of
Applied Mathematics, The Hong Kong Polytechnic University, Hung Hom,
Kowloon, Hong Kong. This author's work was partially supported by the National Natural Science Foundation of China (Grant
No. 11101303).},\hspace{4mm} Zheng-Hai Huang \thanks{Email:
huangzhenghai@tju.edu.cn. Department of Mathematics, School of
Science, Tianjin University, Tianjin, China. This author's work was
supported by the National Natural Science Foundation of China (Grant No. 11171252).}, \hspace{4mm} Chen Ling \thanks{Email:
macling@hdu.edu.cn. School of Science, Hangzhou Dianzi University, Hangzhou 310018, China. This author's work was
supported by the National Natural Science Foundation of China (Grant No. 10971187
and Grant No. 11171083), and the Zhejiang Provincial Natural Science Foundation of China (Grant No. Y6100366).},
\hspace{4mm} Liqun Qi \thanks{Email:
maqilq@polyu.edu.hk. Department of Applied Mathematics, The Hong
Kong Polytechnic University, Hung Hom, Kowloon, Hong Kong. This
author's work was supported by the Hong Kong Research Grant
Council.} }

\date{September 12, 2011} \maketitle

\begin{abstract}
We generalize the concept of the {\em symmetric hyperdeterminants}
for symmetric tensors to the {\em E-determinants} for general
tensors. 
We show that the E-determinant inherits many properties of the
determinant of a matrix. These properties include: solvability of
polynomial systems, the E-determinat of the composition of tensors,
product formula for the E-determinant of a block tensor, Hadamard's
inequality, Ger\v{s}grin's inequality and Minikowski's inequality.
As a simple application, we show that if the leading coefficient
tensor of a polynomial system is a triangular tensor with nonzero
diagonal elements, then the system definitely has a solution.  We
investigate the characteristic polynomial of a tensor through the
E-determinant. Explicit formulae for the coefficients of the
characteristic polynomial are given when the dimension is two.

\noindent {\bf Key words:}\hspace{2mm} Tensor,
eigenvalue, determinant, characteristic polynomial \vspace{3mm}

\end{abstract}

\newpage
\section{Introduction}
\setcounter{Assumption}{0} \setcounter{Theorem}{0}
\setcounter{Proposition}{0} \setcounter{Corollary}{0}
\setcounter{Lemma}{0} \setcounter{Definition}{0}
\setcounter{Remark}{0} \setcounter{Algorithm}{0}
\setcounter{Example}{0} \hspace{4mm} Eigenvalues of tensors,
proposed by Qi \cite{q05} and Lim \cite{l0} independently in 2005,
have attracted much attention in the literature and found various
applications in science and engineering, see
\cite{bv,cs,cpz,cpz1,cglm,hq,lq,l,nqz,q06,q07,q11,qyw,rv,rp,yy,zqy}
and references therein. The concept of {\em symmetric
hyperdeterminant} was introduced by Qi \cite{q05} to investigate the
eigenvalues of a symmetric tensor. Let ${\cal T}=(t_{i_1\ldots
i_m})$ be an $m$-th order $n$-dimensional tensor,
$\mathbf{x}=(x_i)\in\mathbb{C}^n$ (the $n$-dimensional complex
space) and ${\cal T}\mathbf{x}^{m-1}$ be an $n$-dimensional vector
with its $i$-th element being
$\sum_{i_2=1}^n\cdots\sum_{i_m=1}^nt_{ii_2\ldots i_m}x_{i_2}\cdots
x_{i_m}$. Then, when ${\cal T}$ is symmetric, its symmetric
hyperdeterminant is the resultant $\mbox{Res}({\cal
T}\mathbf{x}^{m-1})$ (for the definition of the resultant, see the
next section). The symmetric hyperdeterminant of a symmetric tensor
is equal to the product of all of the eigenvalues of that tensor
\cite{q05}. Recently, Li, Qi and Zhang \cite{lq} proved that the
constant term of the E-characteristic polynomial of tensor ${\cal
T}$ (not necessarily symmetric) is a power of
the resultant $\mbox{Res}({\cal T}\mathbf{x}^{m-1})$. They further
found that the resultant $\mbox{Res}({\cal T}\mathbf{x}^{m-1})$ is
an invariant of ${\cal T}$ under the orthogonal linear
transformation group. Li, Qi and Zhang \cite{lq} pointed out that
the resultant $\mbox{Res}({\cal T}\mathbf{x}^{m-1})$ deserves
further study, since it has close relation to the eigenvalue theory
of tensors. In this paper, we study $\mbox{Res}({\cal
T}\mathbf{x}^{m-1})$ systematically. Note that $\mbox{Res}({\cal
T}\mathbf{x}^{m-1})$ is different from the hyperdeterminant
investigated in \cite{b,c1,c2,clo,d,gkz,gkz1,m,ms,ms1,sl,sy}. We now
give the following definition.

\begin{Definition}\label{def-3}
Let ${\cal T}\in\mathbb{T}(\mathbb{C}^n,m)$ (the space of
$m$-th order $n$-dimensional tensors).  Then its {\em
E-determinant}, denoted by $\mbox{Edet}({\cal T})$, is defined as
the resultant of the polynomial system ${\cal
T}\mathbf{x}^{m-1}=\mathbf{0}$.
\end{Definition}
Here we use the prefix ``E" to highlight the relation of
$\mbox{Edet}({\cal T})$ with the eigenvalue theory of tensors.

The rest of this paper is organized as follows.

In the next section, we present some basic properties of the E-determinant. Then, in Section 3, we
show that the solvability of a polynomial system is
characterized by the E-determinant of the leading coefficient
tensor of the polynomial system.

A tensor ${\cal T}\in\mathbb{T}(\mathbb{C}^n,m)$ induces a
homogenous polynomial map $\hat {\cal T} : \mathbb{C}^n \to \mathbb{C}^n$  as
$\mathbf{x}\mapsto {\cal T}\mathbf{x}^{m-1}$. Let ${\cal
U}\in\mathbb{T}(\mathbb{C}^n,p)$ and ${\cal
V}\in\mathbb{T}(\mathbb{C}^n,q)$ for $p, q \ge 2$, and homogenous polynomial
maps $\hat {\cal U}$ and $\hat {\cal V}$ be induced by $\cal U$ and
$\cal V$ respectively. In Section 4, we show that the composition of
$\hat {\cal U}$ and $\hat {\cal V}$ is another homogenous polynomial map
$\hat {\cal W}$ induced by a tensor ${\cal
W}\in\mathbb{T}\left(\mathbb{C}^n,(p-1)(q-1)+1\right)$.   We show
that $\mbox{Edet}({\cal W}) = 0$ if and only if  $\mbox{Edet}({\cal
U})\mbox{Edet}({\cal V}) = 0$.   We conjecture that
$$\mbox{Edet}({\cal W}) = \left(\mbox{Edet}({\cal
U})\right)^{(q-1)^{n-1}}\left(\mbox{Edet}({\cal
V})\right)^{(p-1)^n}$$ and prove that this conjecture is true when
$\min \{ p, q \} = 2$.


Block tensors are discussed in Section 5. We give an expression of
the E-determinant of a tensor which has an ``upper triangular
structure", based on the E-determinants of its two diagonal
sub-tensors.

As a simple application of the E-determinant theory, in Section 6,
we show that if the leading coefficient tensor of a polynomial
system is a triangular tensor with nonzero diagonal elements, then
the system definitely has a solution.

Based upon a result of Morozov and Shakirov \cite{ms}, in Section 7,
we give a trace formula for the E-determinant.  This formula
involves some differential operators. Using this formula, we will
establish an explicit formula for the E-determinant when the
dimension is two.   As this needs to use some results in Section 8,
we will do this in Section 9.

The E-determinant contributes to the characteristic polynomial
theory of tensors. In Section 8, we analyze various related
properties of the characteristic polynomial and the E-determinant.
Especially, a trace formula for the characteristic polynomial is
presented, which has potential applications in various areas, such
as scientific computing and geometrical analysis of eigenvalues. We
also generalize the eigenvalue representation for the determinant of
a matrix to the E-determinant of a tensor. Under an assumption
(Assumption \ref{assump}),  we transform the positive
semidefiniteness problem of an even order tensor to a computable
condition (see Proposition \ref{prop-4-1}).

In Section 9, we give explicit formulae for the E-determinant and
the characteristic polynomial when the dimension is two.

We generalize some inequalities of the determinant for a matrix to
the E-determinant for a tensor in Section 10. Among them, we present
generalizations of Hadamard's inequality, Ger\v{s}grin's inequality
and Minikowski's inequality. These inequalities give estimations for
the E-determinant in terms of the entries of the underlying tensor.

Some final remarks are given in Section 11.

The following is the notation that is used in the sequel. Scalars
are written as lowercase letters ($\lambda,a,\ldots$); vectors are
written as  bold lowercase letters ($\mathbf{x}=(x_i)$, $\ldots$);
matrices are written as italic capitals ($A=(a_{ij})$, $\ldots$);
tensors are written as calligraphic letters (${\cal T}=(t_{i_1\ldots
i_m})$, $\ldots$); and, sets are written as blackboard bold letters
($\mathbb{T},\mathbb{S},\ldots$).

Given a ring $\mathbb{K}$ (hereafter, we mean a commutative ring with $1$ \cite[Pages 83-84]{lang}, e.g., $\mathbb{C}$),
we denote by $\mathbb{K}[\mathbb{E}]$ the polynomial ring consists of polynomials in the set of indeterminate $\mathbb{E}$ with coefficients in $\mathbb{K}$. Especially, we denote by $\mathbb{K}[{\cal T}]$ the polynomial ring consists of polynomials in indeterminate $\{t_{i_1\ldots i_m}\}$ with coefficients in $\mathbb{K}$, and similarly for $\mathbb{K}[\lambda]$, $\mathbb{K}[A]$, $\mathbb{K}[\lambda,{\cal T}]$, etc.

For a matrix $A$, $A^T$ denotes its transpose and $\mbox{Tr}(A)$ denotes its trace. We denote by $\mathbb{N}_+$ the set
of all positive integers and $\mathbf{e}_i$ the $i$-th identity vector, i.e., the $i$-th column vector of the identity matrix $E$. Throughout this paper, unless stated otherwise, integers $m,n\geq 2$ and tensors refer to $m$-th order $n$-dimensional tensors with entries in $\mathbb{C}$.

\section{Basic Properties of the E-Determinant}
\setcounter{Assumption}{0} \setcounter{Theorem}{0}
\setcounter{Proposition}{0} \setcounter{Corollary}{0}
\setcounter{Lemma}{0} \setcounter{Definition}{0}
\setcounter{Remark}{0} \setcounter{Algorithm}{0}
\setcounter{Example}{0} \hspace{4mm}

Let ${\cal E}$ be the identity tensor of
appropriate order and dimension, e.g., $e_{i_1\ldots i_m}=1$ if and only if $i_1=\cdots=i_m\in\{1,\ldots,n\}$,
and zero otherwise.
The following definitions were introduced by Qi \cite{q05}.
\begin{Definition}\label{def-1}
Let ${\cal T}\in\mathbb{T}(\mathbb{C}^n,m)$.
For some $\lambda\in\mathbb{C}$, if system $\left(\lambda {\cal E}-{\cal T}\right)\mathbf{x}^{m-1}=\mathbf{0}$ has a
solution $\mathbf{x}\in\mathbb{C}^n\setminus\{\mathbf 0\}$, then $\lambda$ is called an
eigenvalue of tensor ${\cal T}$ and ${\mathbf x}$ an eigenvector of ${\cal T}$ associated with $\lambda$.
\end{Definition}
We denote by $\sigma({\cal T})$ the set of all eigenvalues of tensor ${\cal T}$.

\begin{Definition}\label{def-2}
Let ${\cal T}\in\mathbb{T}(\mathbb{C}^n,m)$. The E-determinant of $\lambda{\cal E}-{\cal T}$ which is a
polynomial in $\left(\mathbb{C}[{\cal T}]\right)[\lambda]$, denoted by $\psi(\lambda)$, is called the {\em characteristic polynomial} of tensor ${\cal T}$.
\end{Definition}
If $\lambda$ is a root of $\psi(\lambda)$ of multiplicity $s$, then we call $s$ the {\em algebraic multiplicity} of eigenvalue $\lambda$.
Denote by ${\mathbb{V}}(f)$ the variety of the principal ideal $\langle f\rangle$ generated by $f$ \cite{clo,clo1,lang}. Then, we have the following result.

\begin{Theorem}\label{thm-1}
Let ${\cal T}\in \mathbb{T}(\mathbb{C}^n,m)$. Then $\psi\in \mathbb{C}[\lambda,{\cal T}]$ is homogenous of degree $n(m-1)^{n-1}$ and
\begin{eqnarray}\label{thm-1-1}
{\mathbb{V}}(\psi(\lambda))=\sigma({\cal T}).
\end{eqnarray}
\end{Theorem}

When ${\cal T}$ is symmetric, Qi proved \reff{thm-1-1} in \cite[Theorem 1(a)]{q05}.

For $f\in\mathbb{K}[\mathbf{x}]$, we denote by $\mbox{deg}(f)$ the
degree of $f$. If every monomial in $f$ has degree $\mbox{deg}(f)$,
then $f$ is called homogenous of degree $\mbox{deg}(f)$. Given a
system of polynomials $h:=\{h_1,\ldots,h_n\}$ with
$h_i\in\mathbb{C}[\mathbf{x}]$ being homogenous of degree
$r_i\in\mathbb{N}_+$. The {\em resultant} of the polynomial system
$h$, denoted by $\mbox{Res}_{r_1,\ldots,r_n}(h)$ or simply
$\mbox{Res}(h)$, is defined as an irreducible polynomial in the
coefficients of $h$ such that $\mbox{Res}(h)=0$ if and only if
$h=\mathbf{0}$ has a solution in
$\mathbb{C}^n\setminus\{\mathbf{0}\}$. Furthermore, $\mbox{Res}(h)$
is homogenous of degree $\Pi_{i\neq j}r_i$ in the coefficients of
$h_j$ for every $j\in\{1,\ldots, n\}$. So, it is homogeneous of
total degree $\sum_{i=1}^n\Pi_{j\neq i}^nr_j$ \cite[Proposition
13.1.1]{gkz}, see also \cite[Page 713]{ms}. These, together with Definition \ref{def-3} and
\cite[Theorem 3.2.3(b)]{clo}, immediately imply the following proposition.

\begin{Proposition}\label{prop-000}
Let ${\cal T}\in\mathbb{T}(\mathbb{C}^n,m)$. Then,
\begin{itemize}
\item [(i)] For every $i\in\{1,\ldots,n\}$, let $\mathbb{K}_i:=\mathbb{C}[\{t_{ji_2\ldots i_m}\;|\;j,i_2,\ldots,i_m=1,\ldots,n,\;j\neq i\}]$. Then $\mbox{Edet}({\cal T})\in\mathbb{K}_i[\{t_{ii_2\ldots i_m}\;|\;i_2,\ldots,i_m=1,\ldots,n\}]$
is homogenous of degree $(m-1)^{n-1}$.
\item [(ii)] $\mbox{Edet}({\cal T})\in\mathbb{C}[{\cal T}]$ is irreducible and homogeneous of degree $n(m-1)^{n-1}$.
\item[(iii)] $\mbox{Edet}({\cal E})=1$.
\end{itemize}
\end{Proposition}

By Proposition \ref{prop-000}, we have the following corollary.
\begin{Corollary}\label{cor-00}
Let ${\cal T}\in\mathbb{T}(\mathbb{C}^n,m)$. If for some $i$,
$t_{ii_2\ldots i_m}=0$ for all $i_2,\ldots,i_m\in\{1,\ldots,n\}$,
then $\mbox{Edet}({\cal T})=0$. In particular, the E-determinant of the zero tensor is zero.
\end{Corollary}

\noindent {\bf Proof.} Let $\mathbb{K}_i:=\mathbb{C}[\{t_{ji_2\ldots
i_m}\;|\;j,i_2,\ldots,i_m=1,\ldots,n,\;j\neq i\}]$.   Then by
Proposition \ref{prop-000} (i) $\mbox{Edet}({\cal T})$ is a
homogenous polynomial in the variable set $\{t_{ii_2\ldots
i_m}\;|\;i_2,\ldots,i_m=1,\ldots,n\}$ with coefficients in the ring
$\mathbb{K}_i$. As $t_{ii_2\ldots i_m}=0$ for all
$i_2,\ldots,i_m\in\{1,\ldots,n\}$ by the assumption,
$\mbox{Edet}({\cal T})=0$ as desired. \ep

By Proposition \ref{prop-000} (ii), we have another corollary as
follows.

\begin{Corollary}\label{cor-001}
Let ${\cal T}\in\mathbb{T}(\mathbb{C}^n,m)$ and $\alpha \in
\mathbb{C}$. Then
$$\mbox{Edet}(\alpha{\cal T})= \left(\alpha\right)^{n(m-1)^{n-1}}\mbox{Edet}({\cal
T}).$$
\end{Corollary}

\section{Solvability of Polynomial Equations}
\setcounter{Assumption}{0}
\setcounter{Theorem}{0} \setcounter{Proposition}{0}
\setcounter{Corollary}{0} \setcounter{Lemma}{0}
\setcounter{Definition}{0} \setcounter{Remark}{0}
\setcounter{Algorithm}{0}  \setcounter{Example}{0}
\hspace{4mm} Let matrix $A\in\mathbb{T}(\mathbb{C}^n,2)$, we know that
\begin{itemize}
\item $\mbox{Det}(A)=0$ if and only if $A\mathbf{x}=\mathbf{0}$ has a solution in $\mathbb{C}^n\setminus\{\mathbf{0}\}$; and,
\item $\mbox{Det}(A)\neq 0$ if and only if $A\mathbf{x}=\mathbf{b}$ has a unique solution in $\mathbb{C}^{n}$ for every $\mathbf{b}\in\mathbb{C}^n$.
\end{itemize}
We generalize such a result to the E-determinant and polynomial system in this section.

\begin{Theorem}\label{prop-0}
Let ${\cal T}\in\mathbb{T}(\mathbb{C}^n,m)$. Then,
\begin{itemize}
\item [(i)] $\mbox{Edet}({\cal T})=0$ if and only if ${\cal T}\mathbf{x}^{m-1}=\mathbf{0}$ has a solution in $\mathbb{C}^n\setminus\{\mathbf{0}\}$.
\item [(ii)] $\mbox{Edet}({\cal T})\neq 0$ only if ${\cal T}\mathbf{x}^{m-1}={\cal B}^{m-1}\mathbf{x}^{m-2}+\cdots+{\cal B}^3\mathbf{x}^2+A\mathbf{x}+\mathbf{b}$ has a solution in $\mathbb{C}^n$ for every $\mathbf{b}\in\mathbb{C}^n$, $A\in\mathbb{T}(\mathbb{C}^n,2)$, and ${\cal B}^j\in\mathbb{T}(\mathbb{C}^n,j)$ for $j=3,\ldots,m-1$.
\end{itemize}
\end{Theorem}

\noindent {\bf Proof.} (i) It follows from Definition \ref{def-3} immediately.

(ii) Suppose that $\mbox{Edet}({\cal T})\neq 0$. For any
$\mathbf{b}\in\mathbb{C}^n$, $A\in\mathbb{T}(\mathbb{C}^n,2)$, and
${\cal B}^j\in\mathbb{T}(\mathbb{C}^n,j)$ for $j=3,\ldots,m-1$, we
define tensor ${\cal U}\in\mathbb{T}(\mathbb{C}^{n+1},m)$ as
follows:
\begin{eqnarray}\label{new-1}
u_{i_1i_2\ldots i_m}:=\left\{\begin{array}{cl}t_{i_1i_2\ldots i_m}&\forall i_j\in\{1,\ldots,n\},\;j\in\{1,\ldots,m\},\\
-b_{i_1}&\forall i_1\in\{1,\ldots,n\}\;\mbox{and}\; i_2=\cdots=i_m=n+1,\\-a_{i_1i_2}&\forall i_1,i_2\in\{1,\ldots,n\}\;\mbox{and}\; i_3=\cdots=i_m=n+1,\\-b^k_{i_1\ldots i_k}&\forall i_1,\ldots,i_k\in\{1,\ldots,n\}\;\mbox{and}\; i_{k+1}=\cdots=i_m=n+1,\\&\;\;\;\;\forall k=3,\ldots,m-1,\\0&\mbox{otherwise}.\end{array}\right.
\end{eqnarray}
By Corollary \ref{cor-00}, we have that $\mbox{Edet}({\cal U})=0$.
Hence, by (i), there exists
$\mathbf{y}:=(\mathbf{x}^T,\alpha)^T\in\mathbb{C}^{n+1}\setminus\{\mathbf{0}\}$
such that ${\cal U}\mathbf{y}^{m-1}=\mathbf{0}$. Consequently, by
\reff{new-1} and the first $n$ equations in ${\cal
U}\mathbf{y}^{m-1}=\mathbf{0}$, we know that
\begin{eqnarray}\label{new-2}
{\cal T}\mathbf{x}^{m-1}-\alpha{\cal B}^{m-1}\mathbf{x}^{m-2}-\cdots-\alpha^{m-3}{\cal B}^3\mathbf{x}^2-\alpha^{m-2}A\mathbf{x}-\alpha^{m-1}\mathbf{b}=\mathbf{0}.
\end{eqnarray}
Furthermore, we claim that $\alpha\neq 0$. Otherwise, from \reff{new-2}, ${\cal T}\mathbf{x}^{m-1}=\mathbf{0}$ which means $\mbox{Edet}({\cal T})=0$ by (i).
It is a contradiction.  Hence, from \reff{new-2} we know that $\frac{\mathbf{x}}{\alpha}$ is a solution to
\begin{eqnarray*}
{\cal T}\mathbf{x}^{m-1}={\cal B}^{m-1}\mathbf{x}^{m-2}+\cdots+{\cal B}^3\mathbf{x}^2+A\mathbf{x}+\mathbf{b}.
\end{eqnarray*}
The proof is complete. \ep

So, like the determinants of linear equations, the E-determinants
are criterions for the solvability of non-linear polynomial
equations. It is interesting to investigate whether ${\cal
T}\mathbf{x}^{m-1}={\cal B}^{m-1}\mathbf{x}^{m-2}+\cdots+{\cal
B}^3\mathbf{x}^2+A\mathbf{x}+\mathbf{b}$ has only finitely many
solutions whenever $\mbox{Edet}({\cal T})\neq 0$.

\section{Composition of Homogenous Polynomial Maps}
\setcounter{Assumption}{0} \setcounter{Theorem}{0}
\setcounter{Proposition}{0} \setcounter{Corollary}{0}
\setcounter{Lemma}{0} \setcounter{Definition}{0}
\setcounter{Remark}{0} \setcounter{Algorithm}{0}
\setcounter{Example}{0} \hspace{4mm} Let ${\cal
T}\in\mathbb{T}(\mathbb{C}^n,m)$. Then it induces a homogenous polynomial map
from $\mathbb{C}^n$ to $\mathbb{C}^n$ defined as
$\mathbf{x}\mapsto{\cal T}\mathbf{x}^{m-1}$. We denote this
homogenous polynomial map as $\hat{\cal T}$. Let ${\cal
U}\in\mathbb{T}(\mathbb{C}^n,p)$ and ${\cal
V}\in\mathbb{T}(\mathbb{C}^n,q)$ with $p,q\geq 2$. Then, the
composition map denoted as $\hat {\cal U}\circ \hat {\cal V}$ in the
usual sense is well-defined. Actually, $\hat {\cal U}\circ \hat
{\cal V}:\mathbb{C}^n\ra\mathbb{C}^n$ and $\left(\hat {\cal U}\circ
\hat {\cal V}\right)(\mathbf{x}):=\hat {\cal U}\left(\hat {\cal
V}(\mathbf{x})\right)={\cal U}\left({\cal
V}\mathbf{x}^{q-1}\right)^{p-1}$. Then, it is easy to see that
\begin{eqnarray}\label{in-1}
\left[\left(\hat {\cal U}\circ \hat {\cal V}\right) (\mathbf{x})\right]_i:=\sum_{i_2,\ldots,i_p=1}^nu_{ii_2\ldots i_p}\left({\cal V}\mathbf{x}^{q-1}\right)_{i_2}\cdots \left({\cal V}\mathbf{x}^{q-1}\right)_{i_p},\;\;\forall \mathbf{x}\in\mathbb{C}^n,\;\forall i=1,\ldots,n.
\end{eqnarray}
By the definition of $\hat{\cal T}$ for tensor ${\cal T}$, we can see that there exists a tensor ${\cal W}$ such that $\hat{\cal W}=\hat {\cal U}\circ \hat {\cal V}$.
Furthermore, it is easy to get the following result.

\begin{Proposition}\label{prop-0000}
Let ${\cal U}\in\mathbb{T}(\mathbb{C}^n,p)$ and ${\cal V}\in\mathbb{T}(\mathbb{C}^n,q)$ with $p,q\geq 2$. Then, with ${\cal W}\in\mathbb{T}(\mathbb{C}^n,1+(p-1)(q-1))$ defined as
\begin{eqnarray}\label{new-11}
w_{ij_{(i_2,2)}\ldots j_{(i_2,q)}\ldots j_{(i_p,2)}\ldots j_{(i_p,q)}}:=\sum_{i_2,\ldots,i_p=1}^nu_{ii_2\ldots i_p}v_{i_2j_{(i_2,2)}\ldots j_{(i_2,q)}}\cdots v_{i_pj_{(i_p,2)}\ldots j_{(i_p,q)}},
\end{eqnarray}
we have $\hat{\cal W}=\hat {\cal U}\circ \hat {\cal V}$. So, it is reasonable to define the composition of ${\cal U}$ and ${\cal V}$ as ${\cal U}\circ {\cal V}:={\cal W}$.
\end{Proposition}

Note that when $p=q=2$, i.e., both ${\cal U}$ and ${\cal V}$ are
matrices. Then the composition ${\cal U}\circ {\cal V}$ reduces to
the usual multiplication of matrices, and moreover it is uniquely
defined. Nonetheless, when $m>2$, there exist many ${\cal W}$'s
satisfying $\hat{\cal W}=\hat {\cal U}\circ \hat {\cal V}$. Hence,
we need \reff{new-11} to make the composition of ${\cal U}$ and
${\cal V}$ uniquely defined.

\begin{Proposition}\label{prop-00}
Let ${\cal U}\in\mathbb{T}(\mathbb{C}^n,p)$ and ${\cal V}\in\mathbb{T}(\mathbb{C}^n,q)$ with $p,q\geq 2$. Then, we have
\begin{eqnarray*}
\mbox{Edet}({\cal U}\circ {\cal V})=0 \quad \Longleftrightarrow\quad \mbox{Edet}({\cal U})\mbox{Edet}({\cal V})=0.
\end{eqnarray*}
\end{Proposition}

\noindent {\bf Proof.} The proof for $``\Leftarrow"$: Suppose that
$\mbox{Edet}({\cal V})=0$, by Theorem \ref{prop-0} (i), there exists
$\mathbf{x}\in\mathbb{C}^n\setminus\{\mathbf{0}\}$ such that ${\cal
V}\mathbf{x}^{q-1}=\mathbf{0}$. From \reff{in-1}, it is easy to see
that $\left({\cal U}\circ {\cal V}\right)
\mathbf{x}^{(p-1)(q-1)}=\mathbf{0}$. By Theorem \ref{prop-0} (i)
again, we know that $\mbox{Edet}({\cal U}\circ{\cal V})=0$.

Now, if $\mbox{Edet}({\cal V})\neq 0$, then $\mbox{Edet}({\cal
U})=0$ by the hypothesis, which implies that there exists
$\mathbf{x}\in\mathbb{C}^n\setminus\{\mathbf{0}\}$ such that ${\cal
U}\mathbf{x}^{p-1}=\mathbf{0}$. Since $\mbox{Edet}({\cal V})\neq 0$,
by Theorem \ref{prop-0} (ii), we know that for $\mathbf{x}$ there
exists $\mathbf{y}\in\mathbb{C}^{n}\setminus \{\mathbf{0}\}$ such
that ${\cal V}\mathbf{y}^{q-1}=\mathbf{x}$. Consequently, by
\reff{in-1}, we get that $\left({\cal U}\circ {\cal V}\right)
\mathbf{y}^{(p-1)(q-1)}=\mathbf{0}$. By Theorem \ref{prop-0} (i),
$\mbox{Edet}({\cal U}\circ {\cal V})=0$.

The proof for $``\Rightarrow"$: Suppose that $\mbox{Edet}({\cal
U}\circ {\cal V})=0$. By Theorem \ref{prop-0} (i), there exists
$\mathbf{y}\in\mathbb{C}^n\setminus\{\mathbf{0}\}$ such that
$\left({\cal U}\circ {\cal V}\right)
\mathbf{y}^{(p-1)(q-1)}=\mathbf{0}$. If $\mbox{Edet}({\cal V})=0$,
then we are done. If $\mbox{Edet}({\cal V})\neq 0$, then
$\mathbf{x}:={\cal V}\mathbf{y}^{q-1}\neq \mathbf{0}$. Consequently,
it holds that ${\cal U}\mathbf{x}^{p-1}=\left({\cal U}\circ {\cal
V}\right) \mathbf{y}^{(p-1)(q-1)}=\mathbf{0}$. By Theorem
\ref{prop-0} (i) again, we know that $\mbox{Edet}({\cal U})=0$.
Therefore, $\mbox{Edet}({\cal U})\mbox{Edet}({\cal V})=0$.  \ep

The following corollary is a direct consequence of \cite[Theorem
1(b)]{q05} and Proposition \ref{prop-00}.
\begin{Corollary}\label{cor-0}
Let ${\cal U}\in\mathbb{T}(\mathbb{C}^n,p)$ and ${\cal
V}\in\mathbb{T}(\mathbb{C}^n,q)$ with $p,q\geq 2$. Then, ${\cal
U}\circ {\cal V}$ has zero as its eigenvalue if and only if one of
${\cal U}$ and ${\cal V}$ has zero as its eigenvalue.
\end{Corollary}

We have the following conjecture.

\begin{Conjecture}\label{conj-00}
Let ${\cal U}\in\mathbb{T}(\mathbb{C}^n,p)$ and ${\cal
V}\in\mathbb{T}(\mathbb{C}^n,q)$ with $p,q\geq 2$. Then, we have
\begin{eqnarray*}
\mbox{Edet}({\cal U}\circ {\cal V})=\left(\mbox{Edet}({\cal
U})\right)^{(q-1)^{n-1}}\left(\mbox{Edet}({\cal
V})\right)^{(p-1)^n}.
\end{eqnarray*}
\end{Conjecture}

When $p=q=2$, this conjecture is true as it reduces to the
Cauchy-Binet formula for matrices, i.e, for
$A,B\in\mathbb{T}(\mathbb{C}^n,2)$, it holds that
$\mbox{Det}(AB)=\mbox{Det}(A)\mbox{Det}(B)$ \cite[Page 22]{hj1}. In
the remainder of this section, we show that this conjecture is true
when $\min \{ p, q \} = 2$.

Given a set $\mathbb{E}\subseteq \mathbb{C}^n$, we denote by
$\mathbb{I}(\mathbb{E})\subseteq \mathbb{C}[\mathbf{x}]$ the ideal
of polynomials in $\mathbb{C}[\mathbf{x}]$ which vanishes
identically on $\mathbb{E}$. Given a set of polynomials
$\mathbb{F}:=\{f_1,\ldots,f_s:f_i\in\mathbb{C}[\mathbf{x}]\}$, we
denote by $\mathbb{V}(\mathbb{F})\subseteq \mathbb{C}^n$ the variety
of $\mathbb{F}$, i.e., the set of the common roots of polynomials in
$\mathbb{F}$ \cite{clo,lang}.

The following proposition follows from \reff{new-11}, Proposition \ref{prop-00} and
\cite[Theorem 3.3.5(a)]{clo}.
\begin{Proposition}\label{prop-2-2}
Conjecture \ref{conj-00} is true if $p=2$, i.e., if $G
\in\mathbb{T}(\mathbb{C}^n,2)$ and ${\cal
V}\in\mathbb{T}(\mathbb{C}^n,q)$ for $q\ge 2$, then,
\begin{eqnarray*}
\mbox{Edet}(G \circ {\cal
V})=\left(\mbox{Det}(G)\right)^{(q-1)^{n-1}}\mbox{Edet}({\cal V}).
\end{eqnarray*}
\end{Proposition}

Now, we prove the following result.
\begin{Proposition}\label{prop-2-3} Conjecture \ref{conj-00} is true if $q=2$, i.e., if ${\cal
U}\in\mathbb{T}(\mathbb{C}^n,p)$ and $G
\in\mathbb{T}(\mathbb{C}^n,2)$ for $p\ge 2$, then
\begin{eqnarray*}
\mbox{Edet}({\cal U}\circ G)=\mbox{Edet}({\cal
U})\left(\mbox{Det}(G)\right)^{(p-1)^{n}}.
\end{eqnarray*}
\end{Proposition}

\noindent {\bf Proof.} If $G$ is singular, then the conclusion
follows from Proposition \ref{prop-00}.   We now assume that $G$ is
nonsingular.  By the fact that $\left({\cal U}\circ
G\right)\mathbf{x}^{p-1}={\cal U}\left(G\mathbf{x}\right)^{p-1}$ for
any $\mathbf{x}\in\mathbb{C}^n$ and
$G\in\mathbb{T}(\mathbb{C}^n,2)$, we see that
\begin{eqnarray*}
\mathbb{I}\left(\mathbb{V}(\mbox{Edet}({\cal U}\circ G))\right)=
\mathbb{I}\left(\mathbb{V}(\mbox{Edet}({\cal U}))\right)
\end{eqnarray*}
for any nonsingular $G\in\mathbb{T}(\mathbb{C}^n,2)$. Hence,
$\mbox{Edet}({\cal U}\circ
G)\in\mathbb{I}\left(\mathbb{V}(\mbox{Edet}({\cal U}\circ
G))\right)= \mathbb{I}\left(\mathbb{V}(\mbox{Edet}({\cal
U}))\right)$. Since $\mbox{Edet}({\cal U})\in\mathbb{C}[{\cal U}]$
is irreducible, by Hilbert's Nullstellensatz \cite[Theorem
4.2]{clo1}, we have that
\begin{eqnarray}\label{sec-2-4}
\mbox{Edet}({\cal U}\circ G)=p({\cal U},G)\mbox{Edet}({\cal U})
\end{eqnarray}
for some $p({\cal U},G)\in\mathbb{C}[{\cal U},G]$. Let ${\cal
R}:={\cal U}\circ G$. By Proposition \ref{prop-000} (ii),
$\mbox{Edet}({\cal R})=\mbox{Edet}({\cal U}\circ
G)\in\mathbb{C}[{\cal R}]$ is homogenous of degree $n(p-1)^{n-1}$ in
variables $r_{i_1\ldots i_m}$'s. Hence, by \reff{new-11}, it is
homogenous of degree $n(p-1)^{n-1}$ in variables $u_{i_1\ldots
i_m}$'s and homogenous of degree $n(p-1)^n$ in variables $g_{ij}$'s.
Since $\mbox{Edet}({\cal U})\in\mathbb{C}[{\cal U}]$ is homogenous
of degree $n(p-1)^{n-1}$ by Proposition \ref{prop-000} (ii) again,
$p({\cal U},G)$ is independent of ${\cal U}$ and homogeneous of
degree $n(p-1)^n$ in variables $g_{ij}$'s. Let ${\cal U}={\cal E}$,
it holds that
\begin{eqnarray*}
\left({\cal E}\circ G\right)\mathbf{x}^{p-1}={\cal
E}\left(G\mathbf{x}\right)^{p-1}=\left(\begin{array}{c}(\sum_{j=1}^ng_{1j}x_j)^{p-1}\\\vdots\\(\sum_{j=1}^ng_{nj}x_j)^{p-1}\end{array}\right).
\end{eqnarray*}
Consequently, by \cite[Theorem 3.3.2(b)]{clo}, we have
\begin{eqnarray*}
\mbox{Edet}({\cal E}\circ
G)=\left[\mbox{Res}\left(\begin{array}{c}(\sum_{j=1}^ng_{1j}x_j)\\\vdots\\(\sum_{j=1}^ng_{nj}x_j)\end{array}\right)\right]^{(p-1)^n}=\left(\mbox{Det}(G)\right)^{(p-1)^n}.
\end{eqnarray*}
Therefore, by Proposition \ref{prop-000} (iii) and \reff{sec-2-4},
we conclude that $p({\cal
U},G)=\left(\mbox{Det}(G)\right)^{(p-1)^n}$ and complete the proof.
\ep

The following is a direct corollary of Propositions \ref{prop-2-2} and \ref{prop-2-3}.
\begin{Corollary}\label{cor-4}
Let ${\cal T}\in\mathbb{T}(\mathbb{C}^n,m)$ and
$G\in\mathbb{T}(\mathbb{C}^n,2)$. Then,
\begin{eqnarray*}
\mbox{Edet}(G\circ {\cal
T})\left(\mbox{Det}(G)\right)^{(m-2)(m-1)^{n-1}}=\mbox{Edet}({\cal
T}\circ G).
\end{eqnarray*}
\end{Corollary}

\section{Block Tensors}
\setcounter{Assumption}{0} \setcounter{Theorem}{0}
\setcounter{Proposition}{0} \setcounter{Corollary}{0}
\setcounter{Lemma}{0} \setcounter{Definition}{0}
\setcounter{Remark}{0} \setcounter{Algorithm}{0}
\setcounter{Example}{0} \hspace{4mm}  In the context of matrices, if
a square matrix $A$ can be partitioned as
\begin{eqnarray*}
A=\left(\begin{array}{cc}B&\sharp\\0&C\end{array}\right)
\end{eqnarray*}
with square sub-matrices $B$ and $C$, then
$\mbox{Det}(A)=\mbox{Det}(B)\mbox{Det}(C)$. We now generalize this
property to tensors. The following definition is straightforward.

\begin{Definition}\label{def-4}
Let ${\cal T}\in\mathbb{T}(\mathbb{C}^n,m)$ and $1\leq k\leq n$. Tensor ${\cal U}\in\mathbb{T}(\mathbb{C}^k,m)$ is called a {\em sub-tensor} of ${\cal T}$
associated to the index set $\{j_1,\ldots,j_k\}\subseteq\{1,\ldots,n\}$ if and only if
$u_{i_1\ldots i_m}=t_{j_{i_1}\ldots j_{i_m}}$ for all $i_1,\ldots,i_m\in\{1,\ldots,k\}$.
\end{Definition}

\begin{Theorem}\label{prop-sub}
Let ${\cal T}\in\mathbb{T}(\mathbb{C}^n,m)$. Suppose that there
exists an integer $k$ satisfying $1\leq k\leq n-1$ and $t_{ii_2\ldots
i_m}=0$ for every $i\in\{k+1,\ldots,n\}$ and all indices
$i_2,\ldots,i_m$ such that
$\{i_2,\ldots,i_m\}\cap\{1,\ldots,k\}\neq\emptyset$. Denote by
${\cal U}\in\mathbb{T}(\mathbb{C}^k,m)$ and ${\cal
V}\in\mathbb{T}(\mathbb{C}^{n-k},m)$ the {\em sub-tensors} of ${\cal
T}$ associated to $\{1,\ldots,k\}$ and $\{k+1,\ldots,n\}$,
respectively. Then, it holds that
\begin{eqnarray}\label{sub-1}
\mbox{Edet}({\cal T})=\left[\mbox{Edet}({\cal
U})\right]^{(m-1)^{n-k}}\left[\mbox{Edet}({\cal
V})\right]^{(m-1)^k}.
\end{eqnarray}
\end{Theorem}

\noindent {\bf Proof.}   We first show that
\begin{eqnarray}\label{sub-11}
\mbox{Edet}({\cal T})=0\quad\Longleftrightarrow\quad\mbox{Edet}({\cal
U})\mbox{Edet}({\cal V})=0.
\end{eqnarray}

Suppose that $\mbox{Edet}({\cal T})=0$. Then there exists
$\mathbf{x}\in\mathbb{C}^n\setminus\{\mathbf{0}\}$ such that ${\cal
T}\mathbf{x}^{m-1}=\mathbf{0}$. Denote by
$\mathbf{u}\in\mathbb{C}^k$ the vector consists of $x_1,\ldots,x_k$,
and $\mathbf{v}\in\mathbb{C}^{n-k}$ the vector consists of
$x_{k+1},\ldots,x_n$. If $\mathbf{v}\neq\mathbf{0}$, then
$\mbox{Edet}({\cal V})=0$, and otherwise $\mbox{Edet}({\cal U})=0$.
Hence, we have
\begin{eqnarray*}
\mbox{Edet}({\cal T})=0\quad\Longrightarrow\quad \mbox{Edet}({\cal
U})\mbox{Edet}({\cal V})=0.
\end{eqnarray*}

Conversely, suppose that $\mbox{Edet}({\cal U})\mbox{Edet}({\cal
V})=0$. If $\mbox{Edet}({\cal U})=0$, then there exists
$\mathbf{u}\in\mathbb{C}^k\setminus \{\mathbf{0}\}$ such that ${\cal
U}\mathbf{u}^{m-1}=\mathbf{0}$. Denote
$\mathbf{x}:=(\mathbf{u}^T,\mathbf{0})^T\in\mathbb{C}^n\setminus\{\mathbf{0}\}$,
then ${\cal T}\mathbf{x}^{m-1}=\mathbf{0}$, which implies
$\mbox{Edet}({\cal T})=0$ by Theorem \ref{prop-0} (i). If
$\mbox{Edet}({\cal U})\neq 0$, then $\mbox{Edet}({\cal V})=0$, which
implies that there exists
$\mathbf{v}\in\mathbb{C}^{n-k}\setminus\{\mathbf{0}\}$ such that
${\cal V}\mathbf{v}^{m-1}=\mathbf{0}$. Now, by vector $\mathbf{v}$
and tensor ${\cal T}$, we construct vector
$\mathbf{b}\in\mathbb{C}^n$ as
\begin{eqnarray}\label{new-8}
b_i:=\sum_{j_2,\ldots,j_m=k+1}^nt_{ij_2\ldots j_m}v_{j_2-k}\cdots
v_{j_m-k},\;\forall i\in\{1,\ldots,n\};
\end{eqnarray}
matrix $A\in\mathbb{T}(\mathbb{C}^n,2)$ as
\begin{eqnarray}\label{new-9}
a_{ij}:=\sum_{(q_2,\ldots,q_m)\in\mathbb{D}}t_{iq_2\ldots
q_m}v_{q_2-k}\cdots v_{q_m-k},\;\forall i,j\in\{1,\ldots,n\}
\end{eqnarray}
with $\mathbb{D}:=\{(q_2,\ldots,q_m)\;|\;q_p=j,\;\mbox{for
some}\;p=2,\ldots,m,\;\mbox{and}\;q_l=k+1,\ldots,n,l\neq p\}$; and,
tensors ${\cal B}^s\in\mathbb{T}(\mathbb{C}^n,s)$ for
$s=3,\ldots,m-1$ as
\begin{eqnarray}\label{new-10}
b^s_{ij_2\ldots
j_s}:=\sum_{(q_2,\ldots,q_m)\in\mathbb{D}^s}t_{iq_2\ldots
q_m}v_{q_2-k}\cdots v_{q_m-k},\;\forall
i,j_2,\ldots,j_m\in\{1,\ldots,n\}
\end{eqnarray}
with
\begin{eqnarray*}
\mathbb{D}^s&:=&\{(q_2,\ldots,q_m)\;|\;\{q_{t_2},\ldots, q_{t_s}\}=\{j_2,\ldots,j_s\} \;\mbox{for some pairwise different}\; t_2,\ldots,t_s\\
&&\;\;\;\;\;\;\;\;\;\;\mbox{in}\;\{2,\ldots,m\},\;\mbox{and}\;q_l=k+1,\ldots,n,l\notin\{t_2,\ldots,t_s\}\}.
\end{eqnarray*}

Since $\mbox{Edet}({\cal U})\neq 0$, by Theorem \ref{prop-0}
(ii),
\begin{eqnarray*}
{\cal U}\mathbf{u}^{m-1}+{\cal B}^{m-1}\mathbf{u}^{m-2}+\cdots+{\cal
B}^3\mathbf{u}^2+A\mathbf{u}+\mathbf{b}=\mathbf{0}
\end{eqnarray*}
has a solution $\mathbf{u}\in\mathbb{C}^k$. Let
$\mathbf{x}:=(\mathbf{u}^T,\mathbf{v}^T)^T\in\mathbb{C}^n\setminus\{\mathbf{0}\}$
as $\mathbf{v}\in\mathbb{C}^{n-k}\setminus\{\mathbf{0}\}$. By
\reff{new-8}, \reff{new-9} and \reff{new-10}, we have that
\begin{eqnarray*}
\left({\cal T}\mathbf{x}^{m-1}\right)_i=\left( {\cal
U}\mathbf{u}^{m-1}+{\cal B}^{m-1}\mathbf{u}^{m-2}+\cdots+{\cal
B}^3\mathbf{u}^2+A\mathbf{u}+\mathbf{b}\right)_i=0,\;\forall
i=1,\ldots,k.
\end{eqnarray*}
Furthermore,
\begin{eqnarray*}
\left({\cal T}\mathbf{x}^{m-1}\right)_i=\left( {\cal
V}\mathbf{v}^{m-1}\right)_i=0,\;\forall i=k+1,\ldots,n.
\end{eqnarray*}
Consequently, ${\cal T}\mathbf{x}^{m-1}=\mathbf{0}$ which implies
$\mbox{Edet}({\cal T})=0$ by Theorem \ref{prop-0} (i).

Hence, we proved (\ref{sub-11}).
In the following, we show that \reff{sub-1} holds. Since $\mbox{Edet}({\cal U}),
\mbox{Edet}({\cal V})\in\mathbb{C}[{\cal T}]$, by \reff{sub-11}, we have
\begin{eqnarray*}
\mathbb{V}(\mbox{Edet}({\cal U})\mbox{Edet}({\cal V}))=
\mathbb{V}(\mbox{Edet}({\cal T})),
\end{eqnarray*}
which implies that
\begin{eqnarray*}
\mathbb{I}(\mathbb{V}(\mbox{Edet}({\cal
T})))=\mathbb{I}(\mathbb{V}(\mbox{Edet}({\cal U})\mbox{Edet}({\cal
V}))).
\end{eqnarray*}
By Proposition \ref{prop-000} (ii), $\mbox{Edet}({\cal T})$ is
irreducible. This, together with $\mbox{Edet}({\cal
U})\mbox{Edet}({\cal V})\in\mathbb{I}(\mathbb{V}(\mbox{Edet}({\cal
U})\mbox{Edet}({\cal V})))=\mathbb{I}(\mathbb{V}(\mbox{Edet}({\cal
T})))=\langle\mbox{Edet}({\cal T})\rangle$ and Hilbert's
Nullstellensatz \cite[Theorem 4.2]{clo1}, implies
\begin{eqnarray}\label{new-13}
\left(\mbox{Edet}({\cal U})\mbox{Edet}({\cal V})\right)^r=p({\cal
T})\mbox{Edet}({\cal T})
\end{eqnarray}
for some $r\in\mathbb{N}_+$ and $p({\cal T})\in\mathbb{C}[{\cal
T}]$.  So, $\mbox{Edet}({\cal V})$ is a multiplier of $p({\cal
T})\mbox{Edet}({\cal T})$. Since $\mbox{Edet}({\cal
V})\in\mathbb{C}[{\cal V}]$ is irreducible and $\mbox{Edet}({\cal
T})\in\mathbb{K}[{\cal V}]$ is homogenous with
$\mathbb{K}:=\mathbb{C}[\{t_{ii_2\ldots
i_m}\;|\;i=1,\ldots,k,\;i_2,\ldots,i_m=1,\ldots,n\}]$,
$\mbox{Edet}({\cal V})$ is a multiplier of $\mbox{Edet}({\cal T})$.
Similarly, we can show that $\mbox{Edet}({\cal U})$ is a multiplier
of $\mbox{Edet}({\cal T})$. Consequently, the irreducibilities of
$\mbox{Edet}({\cal U})\in\mathbb{C}[{\cal U}]$ and
$\mbox{Edet}({\cal V})\in\mathbb{C}[{\cal V}]$, together with
\reff{new-13}, imply that
\begin{eqnarray*}
\mbox{Edet}({\cal T})=\left(\mbox{Edet}({\cal
U})\right)^{r_1}\left(\mbox{Edet}({\cal V})\right)^{r_2}
\end{eqnarray*}
for some $r_1,r_2\in\mathbb{N}_+$. Comparing the degrees of the both
sides with Proposition \ref{prop-000} (ii), we get \reff{sub-1}. The
proof is complete. \ep

\section{A Simple Application: Triangular Tensors}
\setcounter{Assumption}{0} \setcounter{Theorem}{0}
\setcounter{Proposition}{0} \setcounter{Corollary}{0}
\setcounter{Lemma}{0} \setcounter{Definition}{0}
\setcounter{Remark}{0} \setcounter{Algorithm}{0}
\setcounter{Example}{0} \hspace{4mm}

Let  ${\cal T} = (t_{i_1\cdots i_m})\in\mathbb{T}(\mathbb{C}^n,m)$.
Suppose that $t_{i_1\cdots i_m} \equiv 0$ if $\min \{ i_2, \cdots,
i_m \} < i_1$.   Then ${\cal T}$ is called an upper triangular
tensor.   Suppose that $t_{i_1\cdots i_m} \equiv 0$ if $\max \{ i_2,
\cdots, i_m \} > i_1$.   Then ${\cal T}$ is called a lower
triangular tensor.  If ${\cal T}$ is either upper or lower
triangular, then ${\cal T}$ is called a triangular tensor.  In
particular, a diagonal tensor is a triangular tensor.

By Theorem \ref{prop-sub}, we have the following proposition.

\begin{Proposition}\label{cor-002}
Suppose that ${\cal T} \in\mathbb{T}(\mathbb{C}^n,m)$
is a triangular tensor.  Then
$$\mbox{Edet}({\cal T})= \displaystyle\prod_{i=1}^n (t_{i\ldots i})^{(m-1)^{n-1}}.$$
\end{Proposition}

By Definition \ref{def-2} and the above proposition, we have the
following corollary.

\begin{Corollary}\label{cor-003}
Suppose that ${\cal T}\in\mathbb{T}(\mathbb{C}^n,m)$ is a triangular
tensor. Then
\begin{eqnarray*}
\sigma({\cal T})= \{t_{i\ldots i}\;|\;i=1,\ldots,n\},
\end{eqnarray*}
and the algebraic multiplicity of $t_{i\ldots i}$ is $(m-1)^{n-1}$
for all $i=1,\ldots,n$.
\end{Corollary}

With Theorem \ref{prop-0}, we have the following simple application
of the E-determinant theory.

\begin{Theorem}
Suppose that $\cal T$ is a triangular tensor with nonzero diagonal
elements.  Then ${\cal T}\mathbf{x}^{m-1}={\cal
B}^{m-1}\mathbf{x}^{m-2}+\cdots+{\cal
B}^3\mathbf{x}^2+A\mathbf{x}+\mathbf{b}$ has a solution in
$\mathbb{C}^n$ for every $\mathbf{b}\in\mathbb{C}^n$,
$A\in\mathbb{T}(\mathbb{C}^n,2)$, and ${\cal
B}^j\in\mathbb{T}(\mathbb{C}^n,j)$ for $j=3,\ldots,m-1$.
\end{Theorem}

We can show that the composition of two upper (lower) triangular
tensors is still an upper (lower) triangular tensor, and Conjecture
\ref{conj-00} is true for two upper (lower) triangular tensors. We
do not go to the details.

\section{A Trace Formula of the E-Determinant}
\setcounter{Assumption}{0}
\setcounter{Theorem}{0} \setcounter{Proposition}{0}
\setcounter{Corollary}{0} \setcounter{Lemma}{0}
\setcounter{Definition}{0} \setcounter{Remark}{0}
\setcounter{Algorithm}{0}  \setcounter{Example}{0} \hspace{4mm} For
the determinant, 
one generalization of Newton's identities:
\begin{eqnarray*}
\mbox{Det}(E+A)=\sum_{k=0}^{\infty}\frac{1}{k!}\left(-\sum_{j=1}^{\infty}\frac{(-1)^j}{j}\mbox{Tr}(A^j)\right)^k
\end{eqnarray*}
is of irreplaceable importance in the theory of the determinant.
Very recently, Morozov and Shakirov \cite{ms} generalized it to the
context of the resultant of a homogenous polynomial system. In this
section, we present a trace formula for the E-determinant based on
the result of Morozov and Shakirov \cite{ms}.

Let ${\cal T}\in\mathbb{T}(\mathbb{C}^n,m)$.
Define the following differential operators:
\begin{eqnarray}\label{cha-1}
\hat g_i:=\sum_{i_2=1}^n\cdots\sum_{i_m=1}^nt_{ii_2\ldots i_m}\frac{\partial}{\partial a_{ii_2}}\cdots \frac{\partial}{\partial a_{ii_m}},\;\forall i=1,\ldots,n,
\end{eqnarray}
where $A$ is an auxiliary $n\times n$ variable matrix consists of
elements $a_{ij}$'s. It is clear that for every $i$, $\hat g_i$ is a
differential operator which belongs to the operator algebra
$\mathbb{C}[\partial A]$, here $\partial A$ is the $n\times n$
matrix with elements $\frac{\partial}{\partial a_{ij}}$'s. In order
to make the operations in \reff{cha-1} convenient to use and the
resulting formulae tidy, we reformulate $\hat g_i$ in the following
way:
\begin{eqnarray}\label{cha-1-1}
\hat g_i=:\sum_{1\leq i_2\leq i_3\leq\cdots\leq i_m\leq n}w_{ii_2\ldots i_m}\frac{\partial}{\partial a_{ii_2}}\cdots \frac{\partial}{\partial a_{ii_m}},\;\forall i=1,\ldots,n.
\end{eqnarray}
While, the corresponding $g_i$ is defined as
\begin{eqnarray}\label{cha-1-2}
g_i(\mathbf{x}):=\sum_{1\leq i_2\leq i_3\leq\cdots\leq i_m\leq n}w_{ii_2\ldots i_m}x_{i_2}\cdots x_{i_m},\;\forall i=1,\ldots,n.
\end{eqnarray}
Then, by direct computation we have
\begin{eqnarray*}
\left.(\hat g_i)^s\left(\sum_{j=1}^na_{ij}x_j\right)^t\right|_{A=0}=\left\{\begin{array}{cl}(s(m-1))!\left[g_i(\mathbf{x})\right]^s&
\mbox{if}\; t=s(m-1),\\0&\mbox{otherwise}.\end{array}\right.
\end{eqnarray*}

We now have the following proposition.
\begin{Proposition}\label{prop-2}
Let ${\cal T}\in \mathbb{T}(\mathbb{C}^n,m)$. Then,
\begin{eqnarray*}
\log\mbox{Edet}({\cal E}-{\cal T})=\left.{\displaystyle\prod_{i=1}^n\left[\sum_{k=0}^{\infty}\frac{m-1}{((m-1)k)!}(\hat g_i)^k\right]\log\mbox{Det}(E-A)}\right|_{A=0}.
\end{eqnarray*}
So,
\begin{eqnarray}\label{cha-3}
\mbox{Edet}({\cal E}-{\cal T})&=& \exp\left(\left.\displaystyle\prod_{i=1}^n\left[\sum_{k=0}^{\infty}\frac{m-1}{((m-1)k)!}(\hat g_i)^k\right]\log\mbox{Det}(E-A)\right|_{A=0}\right)\nonumber\\
&=& \exp\left(\sum_{k_1=0}^{\infty}\cdots\sum_{k_n=0}^{\infty}-\mbox{Tr}_{k_1,\ldots,k_n}({\cal T})\right)
\end{eqnarray}
with the graded components defined as
\begin{eqnarray}\label{cha-4}
\mbox{Tr}_{k_1,\ldots,k_n}({\cal T}):=(m-1)^n\displaystyle\prod_{i=1}^n\left[\frac{(\hat g_i)^{k_i}}{((m-1)k_i)!}\right]\frac{\mbox{Tr}(A^{(m-1)(\sum_{i=1}^nk_i)})}{(m-1)(\sum_{i=1}^nk_i)}
\end{eqnarray}
and $\mbox{Tr}_{0,\ldots,0}({\cal T}):=0$.
\end{Proposition}

\noindent {\bf Proof.} The results follow from \cite[Sections
4-7]{ms}.   We omit the details.\ep

Note that we can derive the expansion of the right hand side of
\reff{cha-3} by using multi-Schur polynomials in terms of
$\mbox{Tr}_{k_1,\ldots,k_n}({\cal T})$'s. Then, a trace formula for
$\mbox{Edet}({\cal T})$ can be derived.

Motivated by \cite{cd,ms}, we define the {\em $d$-th trace} of
tensor ${\cal T}$ as
\begin{eqnarray}\label{cha-5}
\mbox{Tr}_{d}({\cal T})&:=&(m-1)^{n-1}\left[\sum_{\sum_{i=1}^nk_i=d}\displaystyle\prod_{i=1}^n\frac{(\hat g_i)^{k_i}}{((m-1)k_i)!}\right]\mbox{Tr}(A^{(m-1)d})\nonumber\\
&=&d\sum_{\sum_{i=1}^nk_i=d}\mbox{Tr}_{k_1,\ldots,k_n}({\cal T}).
\end{eqnarray}
Then, it follows from \reff{cha-3}, \reff{cha-4} and \reff{cha-5} that
\begin{eqnarray}\label{cha-6}
\mbox{Edet}({\cal E}-{\cal T})=\exp\left(\sum_{d=0}^{\infty}-\frac{\mbox{Tr}_d({\cal T})}{d}\right).
\end{eqnarray}
We remark that \reff{cha-6} is a generalization of the well known
identity
\begin{eqnarray*}
\log \mbox{Det}(E-A)=\mbox{Tr}\left(\log (E-A)\right)
\end{eqnarray*}
for $A\in\mathbb{T}(\mathbb{C}^n,2)$, i.e., a square matrix, by
noticing that $\log (E-A)=-\sum_{k=1}^{\infty}\frac{A^k}{k}$. In order to derive the
expansion of the right hand side in \reff{cha-6}, we need Schur
polynomials which are defined as:
\begin{eqnarray}\label{cha-7}
p_0(t_0)=1,\;\mbox{and}\;p_k(t_1,\ldots,t_k):=\sum_{i=1}^k\sum_{d_j>0,\sum_{j=1}^id_j=k}\frac{\Pi_{j=1}^it_{d_j}}{i!},\;\;\forall k\geq 1,
\end{eqnarray}
where $\{t_0,t_1,\ldots\}$ are variables.
Let $t_0=0$, we obtain the following expansion:
\begin{eqnarray*}
\exp\left(\sum_{k=0}^{\infty}t_k\alpha^k\right)=1+\sum_{k=1}^{\infty}p_k(t_1,\ldots,t_k)\alpha^k.
\end{eqnarray*}
This, together with \reff{cha-6}, implies
\begin{eqnarray}\label{cha-8}
\mbox{Edet}({\cal T})&=&\mbox{Edet}({\cal E}-\left({\cal E}-{\cal T}\right))\nonumber\\
&=&1+\sum_{k=1}^{\infty}p_k\left(-\frac{\mbox{Tr}_1({\cal E}-{\cal T})}{1},\ldots,-\frac{\mbox{Tr}_k({\cal E}-{\cal T})}{k}\right).
\end{eqnarray}

Now, we improve \reff{cha-8} to give an expression of $\mbox{Edet}({\cal T})$ with only finitely many terms. To this end, we need the following proposition.
\begin{Proposition}\label{prop-3}
Let ${\cal T}\in \mathbb{T}(\mathbb{C}^n,m)$. Then, the followings hold:
\begin{itemize}
\item [(i)] for every $d\in\mathbb{N}_+$, $\mbox{Tr}_d({\cal T})\in\mathbb{C}[{\cal T}]$ is homogenous of degree $d$;
\item [(ii)]for every $k\in\mathbb{N}_+$, $p_k\left(-\frac{\mbox{Tr}_1({\cal T})}{1},\ldots,-\frac{\mbox{Tr}_k({\cal T})}{k}\right)\in\mathbb{C}[{\cal T}]$
is homogenous of degree $k$; and,
\item [(iii)]  for any integer $k>n(m-1)^{n-1}$, $p_k\left(-\frac{\mbox{Tr}_1({\cal T})}{1},\ldots,-\frac{\mbox{Tr}_k({\cal T})}{k}\right)\in\mathbb{C}[{\cal T}]$ is zero.
\end{itemize}
\end{Proposition}

\noindent {\bf Proof.} (i) By the formulae of $\hat g_i$'s as in \reff{cha-1}, it is easy to see that
\begin{eqnarray*}
\sum_{\sum_{i=1}^nk_i=d}\displaystyle\prod_{i=1}^n\frac{(\hat g_i)^{k_i}}{((m-1)k_i)!}\in\mathbb{C}[{\cal T},\partial A]
\end{eqnarray*}
is homogeneous, and more explicitly, homogenous of degree $d$ in the variable ${\cal T}$ and homogeneous of degree $(m-1)d$ in the variable $\partial A$. It is also known that
\begin{eqnarray}\label{prop-3-1}
\mbox{Tr}(A^k)=\sum_{i_1=1}^n\cdots\sum_{i_k=1}^na_{i_1i_2}a_{i_2i_3}\cdots a_{i_{k-1}i_k}a_{i_ki_1}\in\mathbb{C}[A]
\end{eqnarray}
is homogeneous of degree $k$. These, together with \reff{cha-5}, implies that $\mbox{Tr}_d({\cal T})\in\mathbb{C}[{\cal T}]$ is homogenous of degree $d$ as desired.

(ii) It follows from (i) and the definitions of Schur polynomials as in \reff{cha-7} directly.

(iii) From Proposition \ref{prop-000} (ii), it is clear that
$\mbox{Edet}({\cal E}-{\cal T})$ is an irreducible polynomial which
is homogenous of degree $n(m-1)^{n-1}$ in the entries of ${\cal
B}:={\cal E}-{\cal T}$. Since the entries of ${\cal B}$ consist of
$1$ and the entries of tensor ${\cal T}$, the highest degree of
$\mbox{Edet}({\cal E}-{\cal T})$ viewed as a polynomial in
$\mathbb{C}[{\cal T}]$ is not greater than $n(m-1)^{n-1}$. This,
together with (ii) which asserts that $p_k(-\mbox{Tr}_1({\cal
T}),\ldots,-\mbox{Tr}_k({\cal T}))\in\mathbb{C}[{\cal T}]$ is
homogenous of degree $k$, implies the result (iii).

The proof is complete. \ep

By Proposition \ref{prop-3} and \reff{cha-8}, we immediately get
\begin{eqnarray}\label{new-12}
\mbox{Edet}({\cal T})&=&\mbox{Edet}({\cal E}-({\cal E}-{\cal T}))\nonumber\\
&=&1+\sum_{k=1}^{n(m-1)^{n-1}}p_k\left(-\frac{\mbox{Tr}_1({\cal E}-{\cal T})}{1},\ldots,-\frac{\mbox{Tr}_k({\cal E}-{\cal T})}{k}\right).
\end{eqnarray}

This is a trace formula for the E-determinant.  It involves the
differential operators $\hat g_i$'s.   In Section 9, we will give an
explicit formula when $n=2$.

\section{The Characteristic Polynomial}
 \setcounter{Assumption}{0}
\setcounter{Theorem}{0} \setcounter{Proposition}{0}
\setcounter{Corollary}{0} \setcounter{Lemma}{0}
\setcounter{Definition}{0} \setcounter{Remark}{0}
\setcounter{Algorithm}{0}  \setcounter{Example}{0} \hspace{4mm} 
By Definition \ref{def-2}, for any ${\cal
T}\in\mathbb{T}(\mathbb{C}^n,m)$, its characteristic polynomial is
$\psi(\lambda)=\mbox{Edet}(\lambda {\cal E}-{\cal T})$.
The characteristic polynomial of a tensor was proposed by Qi in \cite{q05}, and investigated by Cooper and Dutle in \cite{cd} very recently for spectral hypergraph theory. Following up Qi \cite{q05}, Morozov and Shakirov \cite{ms} and Cooper and Dutle \cite{cd}, we discuss some properties of the characteristic polynomial of a tensor related to the E-determinant.

\begin{Theorem}\label{thm-2}
Let ${\cal T}\in \mathbb{T}(\mathbb{C}^n,m)$. Then
\begin{eqnarray*}
\psi(\lambda)&=&\mbox{Edet}(\lambda {\cal E}-{\cal T})\nonumber\\
&=&\lambda^{n(m-1)^{n-1}}+\sum_{k=1}^{n(m-1)^{n-1}}\lambda^{n(m-1)^{n-1}-k}p_k\left(-\frac{\mbox{Tr}_1({\cal T})}{1},\ldots,-\frac{\mbox{Tr}_k({\cal T})}{k}\right)\\
&=&\Pi_{\lambda_i\in\sigma({\cal T})}(\lambda-\lambda_i)^{m_i}\nonumber,
\end{eqnarray*}
where $m_i$ is the algebraic multiplicity of eigenvalue $\lambda_i$.
\end{Theorem}

\noindent {\bf Proof.} The first equality follows from Definition \ref{def-2}, and the last one from Theorem \ref{thm-1}.

By Proposition \ref{prop-3} and \reff{new-12}, we can get
\begin{eqnarray}\label{new-3}
\psi(1)=\mbox{Edet}({\cal E}-{\cal T})
&=&1+\sum_{k=1}^{n(m-1)^{n-1}}p_k\left(-\frac{\mbox{Tr}_1({\cal T})}{1},\ldots,-\frac{\mbox{Tr}_k({\cal T})}{k}\right).
\end{eqnarray}
Consequently, when $\lambda\neq 0$,
\begin{eqnarray*}
\psi(\lambda)&=&\mbox{Edet}(\lambda {\cal E}-{\cal T})\\
&=&\lambda^{n(m-1)^{n-1}}\mbox{Edet}({\cal E}-\frac{{\cal T}}{\lambda})\\
&=&\lambda^{n(m-1)^{n-1}}\left[1+\sum_{k=1}^{n(m-1)^{n-1}}p_k\left(-\frac{\mbox{Tr}_1(\frac{{\cal T}}{\lambda})}{1},\ldots,-\frac{\mbox{Tr}_k(\frac{{\cal T}}{\lambda})}{k}\right)\right]\\
&=&\lambda^{n(m-1)^{n-1}}\left[1+\sum_{k=1}^{n(m-1)^{n-1}}\frac{1}{\lambda^k}p_k\left(-\frac{\mbox{Tr}_1({\cal T})}{1},\ldots,-\frac{\mbox{Tr}_k({\cal T})}{k}\right)\right]\\
&=&\lambda^{n(m-1)^{n-1}}+\sum_{k=1}^{n(m-1)^{n-1}}\lambda^{n(m-1)^{n-1}-k}p_k\left(-\frac{\mbox{Tr}_1({\cal T})}{1},\ldots,-\frac{\mbox{Tr}_k({\cal T})}{k}\right).
\end{eqnarray*}
Here the second equality comes from Corollary \ref{cor-001}; the
third from \reff{new-3}; and, the fourth from Proposition
\ref{prop-3}. Hence, the result follows from the fact that the field
$\mathbb{C}$ is of characteristic zero. The proof is complete.\ep



Theorem \ref{thm-2} gives a trace formula for the characteristic
polynomial of tensor ${\cal T}$ as well as an eigenvalue
representation for it.

For the sequel analysis, we present the following hypothesis.
\begin{Assumption}\label{assump}
Let ${\cal T}\in \mathbb{S}(\mathbb{R}^n,m)$. Suppose that for every
negative eigenvalue $\lambda$ of ${\cal T}$, it possesses a real
eigenvector associated to $\lambda$.
\end{Assumption}
When $m=2$, Assumption \ref{assump} holds. Note that the tensors in \cite[Examples 1 and 2]{q05} satisfy Assumption \ref{assump}.

A tensor ${\cal
T}\in\mathbb{S}(\mathbb{R}^n,m)$ is called positive semidefinite if
and only if $\mathbf{x}^T\left({\cal T}\mathbf{x}^{m-1}\right)\geq
0$ for all $\mathbf{x}\in\mathbb{R}^n$. Obviously, $m$ being even is necessary for positive semidefinite tensors.

\begin{Lemma}\label{lem-8}
Let ${\cal T}\in \mathbb{S}(\mathbb{R}^n,m)$ and $m$ be even. Suppose that Assumption \ref{assump} holds. Then, ${\cal T}$ is positive semidefinite if and only if all the real eigenvalues of ${\cal T}$ are nonnegative.
\end{Lemma}

\noindent {\bf Proof.} If all the real eigenvalues of ${\cal T}$ are nonnegative, then ${\cal T}$ is positive semidefinite by \cite[Theorem 5]{q05}.

If ${\cal T}$ is positive semidefinite and it has a negative eigenvalue, then ${\cal T}$ has a real eigenvector associated to this eigenvalue by the hypothesis. As ${\cal T}$ is symmetric, this eigenvector violates the definition of positive semidefiniteness by \cite[Theorem 3]{q05}. Consequently, if ${\cal T}$ is positive semidefinite, then all its real eigenvalues are nonnegative. \ep

The following classical result on real roots of a polynomial is D\'{e}scartes's Rule of Signs \cite[Theorem 1.5]{s-b}.
\begin{Lemma}\label{lem-7}
The number of positive real roots of a polynomial is at most the number of sign changes in its coefficients.
\end{Lemma}

Let $\mbox{sgn}(\cdot)$ be the sign function for scalars, i.e.,
$\mbox{sgn}(\alpha)=1$ if $\alpha>0$, $\mbox{sgn}(0)=0$ and
$\mbox{sgn}(\alpha)=-1$ if $\alpha<0$.

\begin{Corollary}\label{cor-5}
Let ${\cal T}\in \mathbb{T}(\mathbb{R}^n,m)$, and
\begin{eqnarray*}
\psi(\lambda)&=&\lambda^{n(m-1)^{n-1}}+\sum_{k=1}^{n(m-1)^{n-1}}\lambda^{n(m-1)^{n-1}-k}p_k\left(-\frac{\mbox{Tr}_1({\cal T})}{1},\ldots,-\frac{\mbox{Tr}_k({\cal T})}{k}\right).
\end{eqnarray*}
Suppose that
\begin{eqnarray*}
\mbox{sgn}\left(p_k\left(-\frac{\mbox{Tr}_1({\cal T})}{1},\ldots,-\frac{\mbox{Tr}_k({\cal T})}{k}\right)\right)=(-1)^{k}
\end{eqnarray*}
for all $p_k\left(-\frac{\mbox{Tr}_1({\cal T})}{1},\ldots,-\frac{\mbox{Tr}_k({\cal T})}{k}\right)\neq 0$ with $1\leq k\leq n(m-1)^{n-1}$. Then, all the real roots of $\psi$ are nonnegative.
\end{Corollary}

\noindent {\bf Proof.} Suppose that $m$ is even and $n$ is odd. Then, $n(m-1)^{n-1}$ is odd. Consequently,
\begin{eqnarray*}
\phi(\lambda):=\psi(-\lambda)&=&-\lambda^{n(m-1)^{n-1}}+\sum_{k=1}^{n(m-1)^{n-1}}(-1)^{k+1}\lambda^{n(m-1)^{n-1}-k}p_k\left(-\frac{\mbox{Tr}_1({\cal T})}{1},\ldots,-\frac{\mbox{Tr}_k({\cal T})}{k}\right).
\end{eqnarray*}
Then, by Lemma \ref{lem-7}, $\phi$ defined as above has no positive real root, since the sign of the coefficient of $\phi$ is negative when it is nonzero. Hence, $\psi$ has no negative real root.

The proofs for the other cases for $m$ and $n$ are similar. Consequently, the result follows. The proof is complete. \ep

\begin{Proposition}\label{prop-4-1}
Let $m$ be even, ${\cal T}\in \mathbb{S}(\mathbb{R}^n,m)$, and
\begin{eqnarray*}
\psi(\lambda)&=&\lambda^{n(m-1)^{n-1}}+\sum_{k=1}^{n(m-1)^{n-1}}\lambda^{n(m-1)^{n-1}-k}p_k\left(-\frac{\mbox{Tr}_1({\cal T})}{1},\ldots,-\frac{\mbox{Tr}_k({\cal T})}{k}\right).
\end{eqnarray*}
Suppose that Assumption \ref{assump} holds. Then, ${\cal T}$ is positive semidefinite if
\begin{eqnarray}\label{new-14}
\mbox{sgn}\left(p_k\left(-\frac{\mbox{Tr}_1({\cal T})}{1},\ldots,-\frac{\mbox{Tr}_k({\cal T})}{k}\right)\right)=(-1)^{k}
\end{eqnarray}
for all $p_k\left(-\frac{\mbox{Tr}_1({\cal T})}{1},\ldots,-\frac{\mbox{Tr}_k({\cal T})}{k}\right)\neq 0$ with $1\leq k\leq n(m-1)^{n-1}$.
Furthermore, when $n$ is even and all the complex eigenvalues of ${\cal T}$ have nonnegative real parts, then ${\cal T}$ is positive semidefinite if and only if \reff{new-14} holds for all $p_k\left(-\frac{\mbox{Tr}_1({\cal T})}{1},\ldots,-\frac{\mbox{Tr}_k({\cal T})}{k}\right)\neq 0$ with $1\leq k\leq n(m-1)^{n-1}$.
\end{Proposition}

\noindent {\bf Proof.}
The first result follows from Lemma \ref{lem-8} and Corollary \ref{cor-5} directly.

Now, we prove the second result.
By Lemma \ref{lem-8}, all the real roots of polynomial $\psi$ are nonnegative. Consequently, all the real roots of polynomial $\phi(\lambda):=\psi(-\lambda)$ are nonpositive. Suppose that $\phi$ has negative roots $\{-\alpha_1,\ldots,-\alpha_s\}$ with the corresponding multiplicity set $\{m_1,\ldots,m_s\}$; zero root with multiplicity $k$; and, complex root pairs $\{(\mu_1,\bar \mu_1),\ldots,(\mu_t,\bar \mu_t)\}$ with the corresponding multiplicity set $\{r_1,\ldots,r_s\}$. Consequently,
\begin{eqnarray*}
\phi(\lambda)=\lambda^k\Pi_{i=1}^s(\lambda+\alpha_i)^{m_i}\Pi_{i=1}^t(\lambda^2+(\mu_i+\bar
\mu_i)\lambda+|\mu_i|^2)^{r_i}.
\end{eqnarray*}
Hence, all the coefficients of $\phi$ are nonnegative by the assumption that all the complex eigenvalues of ${\cal T}$ have nonnegative real parts. Moreover, $n(m-1)^{n-1}$ is even, since $n$ is even. These, together with the definition of $\phi$, imply \reff{new-14} immediately.

The proof is complete.\ep

\begin{Remark}\label{rmk-2}
{\em Proposition \ref{prop-4-1} is meaningful: we do not need to
compute all the real eigenvalues of ${\cal T}$ associated with real
eigenvectors, even the smallest, as \cite[Theorem 5]{q05}. All we
need is to check the condition \reff{new-14}, 
when the additional hypotheses in Proposition \ref{prop-4-1} are
satisfied.
}
\end{Remark}


Here are some properties concerning the coefficients of $\psi(\lambda)$.

\begin{Proposition}\label{prop-4}
Let ${\cal T}\in \mathbb{T}(\mathbb{C}^n,m)$. Then,
\begin{itemize}
\item [(i)] $p_1(-\mbox{Tr}_1({\cal T}))=-\mbox{Tr}_1({\cal T})=-(m-1)^{n-1}\sum_{i=1}^nt_{ii\ldots i}$; and,
\item [(ii)] $p_{n(m-1)^{n-1}}\left(-\frac{\mbox{Tr}_1({\cal T})}{1},\ldots,-\frac{\mbox{Tr}_{n(m-1)^{n-1}}({\cal T})}{n(m-1)^{n-1}}\right)=(-1)^{n(m-1)^{n-1}}\mbox{Edet}({\cal T})$.
\end{itemize}
\end{Proposition}

\noindent {\bf Proof.} (i) By \reff{cha-7}, we know that $p_1(-\mbox{Tr}_1({\cal T}))=-\mbox{Tr}_1({\cal T})$. Furthermore, by \reff{cha-5}, it is easy to see that
\begin{eqnarray*}
\mbox{Tr}_1({\cal T})&=&(m-1)^{n-1}\sum_{i=1}^n\frac{\hat g_i}{(m-1)!}\mbox{Tr}(A^{m-1})\\
&=&\frac{(m-1)^{n-1}}{(m-1)!}\sum_{i=1}^n\left[\sum_{i_2=1}^n\cdots\sum_{i_m=1}^nt_{ii_2\ldots i_m}\frac{\partial}{\partial a_{ii_2}}\cdots \frac{\partial}{\partial a_{ii_m}}\right]\mbox{Tr}(A^{m-1})\\
&=&\frac{(m-1)^{n-1}}{(m-1)!}\sum_{i=1}^n\left[\sum_{i_2=1}^n\cdots\sum_{i_m=1}^nt_{ii_2\ldots i_m}\frac{\partial}{\partial a_{ii_2}}\cdots \frac{\partial}{\partial a_{ii_m}}\right]\\
&&\;\;\;\;\cdot\left(\sum_{i_1=1}^n\cdots\sum_{i_{m-1}=1}^na_{i_1i_2}a_{i_2i_3}\cdots a_{i_{m-2}i_{m-1}}a_{i_{m-1}i_1}\right)\\
&=&\frac{(m-1)^{n-1}}{(m-1)!}\sum_{i=1}^n\left[t_{ii\ldots i}\frac{\partial}{\partial a_{ii}}\cdots \frac{\partial}{\partial a_{ii}}(a_{ii})^{m-1}\right]\\
&=&(m-1)^{n-1}\sum_{i=1}^nt_{ii\ldots i}.
\end{eqnarray*}
Here, the fourth equality follows from the fact that: (a) the differential operator in the right hand side of the third equality contains only items $\frac{\partial}{\partial a_{i\star}}$'s for $\star\in\{1,\ldots,n\}$ and the total degree is $m-1$; and, (b) only terms in $\mbox{Tr}(A^{m-1})$ that contain the same $\frac{\partial}{\partial a_{i\star}}$'s of total degree $m-1$ can contribute to the result and this case occurs only when every $\star=i$ by \reff{prop-3-1}.
Consequently, the result (i) follows.

(ii) By Theorem \ref{thm-2}, it is clear that
\begin{eqnarray*}
\psi(0)=\mbox{Edet}(-{\cal T})=p_{n(m-1)^{n-1}}\left(-\frac{\mbox{Tr}_1({\cal T})}{1},\ldots,-\frac{\mbox{Tr}_{n(m-1)^{n-1}}({\cal T})}{n(m-1)^{n-1}}\right).
\end{eqnarray*}
Moreover, $\mbox{Edet}(-{\cal T})\in\mathbb{C}[{\cal T}]$ is homogenous of degree $n(m-1)^{n-1}$ by Proposition \ref{prop-000} (ii), which implies $\mbox{Edet}(-{\cal T})=(-1)^{n(m-1)^{n-1}}\mbox{Edet}({\cal T})$. Consequently, the result follows.   \ep

\begin{Corollary}\label{cor-1}
Let ${\cal T}\in \mathbb{T}(\mathbb{C}^n,m)$. Then,
\begin{itemize}
\item [(i)] $\sum_{\lambda_i\in\sigma({\cal T})}m_i\lambda_i=(m-1)^{n-1}\sum_{i=1}^nt_{ii\ldots i}=\mbox{Tr}_1({\cal T})$, and
\item [(ii)] $\Pi_{\lambda_i\in\sigma({\cal T})}\lambda_i^{m_i}=\mbox{Edet}({\cal T})$.
\end{itemize}
Here $m_i$ is the algebraic multiplicity of eigenvalue $\lambda_i$.
\end{Corollary}

\noindent {\bf Proof.} The results follow from the eigenvalue representation of $\psi(\lambda)$ in Theorem \ref{thm-2} and the coefficients of $\psi(\lambda)$ in Proposition \ref{prop-4} immediately. \ep

\begin{Remark}\label{rmk-1}
{\em In \cite{q05}, Qi proved the results in Corollary \ref{cor-1} for ${\cal T}\in \mathbb{S}(\mathbb{R}^n,m)$. By Theorem \ref{prop-0} and Corollary \ref{cor-1}, we see that the solvability of homogeneous polynomial equations is characterized by the zero eigenvalue of the underlying tensor}.
\end{Remark}

\section{Explicit Formulae When $n=2$}
 \setcounter{Assumption}{0}
\setcounter{Theorem}{0} \setcounter{Proposition}{0}
\setcounter{Corollary}{0} \setcounter{Lemma}{0}
\setcounter{Definition}{0} \setcounter{Remark}{0}
\setcounter{Algorithm}{0}  \setcounter{Example}{0} \hspace{4mm} We
discuss more on the characteristic polynomial $\psi(\lambda)$ and
the E-determinant $\mbox{Edet}({\cal T})$ of a tensor ${\cal
T}\in\mathbb{T}(\mathbb{C}^n,m)$ in this section. Note that the
trace formulae of both the characteristic polynomial and the
E-determinant depend on the $d$-th traces of the underlying tensor
for all $d=1,\ldots,n(m-1)^{n-1}$. Nevertheless, it is very
complicated \cite{cd,ms}. So, we give preliminary results on the
computation of the $d$-th traces of a tensor. In particular, we give
explicit formulae of $\mbox{Tr}_2({\cal T})$ of the tensor $\cal T$ for any order and dimension, and the characteristic
polynomial $\psi(\lambda)$ and the E-determinant $\mbox{Edet}({\cal
T})$ of the tensor $\cal T$ when $n=2$.

The following lemma is a generalization of Proposition \ref{prop-4}
(i).
\begin{Lemma}\label{lem-1}
Let ${\cal T}\in \mathbb{T}(\mathbb{C}^n,m)$. We have
\begin{eqnarray}\label{lem-1-1}
\frac{(\hat g_i)^k}{((m-1)k)!}\mbox{Tr}(A^{(m-1)k})=t_{ii\ldots i}^k
\end{eqnarray}
for all $k\geq 0$ and $i\in\{1,\ldots,n\}$. So,
\begin{eqnarray}\label{lem-1-2}
\sum_{i=1}^n\frac{(\hat
g_i)^k}{((m-1)k)!}\mbox{Tr}(A^{(m-1)k})=\sum_{i=1}^nt_{ii\ldots
i}^k.
\end{eqnarray}
\end{Lemma}

\noindent {\bf Proof.} By \reff{cha-1} and \reff{prop-3-1}, similar
to the proof of Proposition \ref{prop-4-1} (i), we have
\begin{eqnarray*}
(\hat g_i)^k\mbox{Tr}(A^{(m-1)k})
&=&\left[\sum_{i_2=1}^n\cdots\sum_{i_m=1}^nt_{ii_2\ldots i_m}\frac{\partial}{\partial a_{ii_2}}\cdots \frac{\partial}{\partial a_{ii_m}}\right]^k\mbox{Tr}(A^{(m-1)k})\\
&=&\left[\sum_{i_2=1}^n\cdots\sum_{i_m=1}^nt_{ii_2\ldots i_m}\frac{\partial}{\partial a_{ii_2}}\cdots \frac{\partial}{\partial a_{ii_m}}\right]^k\\
&&\;\;\;\;\cdot\left(\sum_{i_1=1}^n\cdots\sum_{i_{(m-1)k}=1}^na_{i_1i_2}a_{i_2i_3}\cdots a_{i_{(m-1)k-1}i_{(m-1)k}}a_{i_{(m-1)k}i_1}\right)\\
&=&\left[t_{ii\ldots i}\frac{\partial}{\partial a_{ii}}\cdots \frac{\partial}{\partial a_{ii}}\right]^k(a_{ii})^{(m-1)k}\\
&=&((m-1)k)!t_{ii\ldots i}^k,
\end{eqnarray*}
which implies \reff{lem-1-1}, and hence \reff{lem-1-2}.  \ep

Before further analysis, we need the following combinatorial result.
\begin{Lemma}\label{lem-2-0}
Let $i\neq j$, $k\geq 1$, $h\geq 1$ and $s\in\{1,\ldots,
\min\{h,k\}(m-1)\}$ be arbitrary but fixed. Then, the number of term
$(a_{ii})^{k(m-1)-s}(a_{ij})^{s}(a_{ji})^{s}(a_{jj})^{h(m-1)-s}$ in
the expansion of $\mbox{Tr}(A^{(k+h)(m-1)})$ is
\begin{eqnarray*}
\left(\begin{array}{c}k(m-1)\\s\end{array}\right)\left(\begin{array}{c}h(m-1)-1\\s-1\end{array}\right)
+\left(\begin{array}{c}h(m-1)\\s\end{array}\right)\left(\begin{array}{c}k(m-1)-1\\s-1\end{array}\right).
\end{eqnarray*}
\end{Lemma}

\noindent {\bf Proof.} For the convenience of the sequel analysis,
we define a {\em packaged element} of $i$ as an ordered collection
of $a_{ij}$, $a_{jj}$'s and $a_{ji}$ with the form:
\begin{eqnarray*}
a_{ij}\underbrace{a_{jj}\cdots a_{jj}}_{p}a_{ji}.
\end{eqnarray*}
The number $p$ of $a_{jj}$'s in a packaged element of $i$ can vary
from $0$ to the maximal number. A packaged element of $j$ can be
defined similarly.

Note that any term in
\begin{eqnarray*}
\mbox{Tr}(A^{(k+h)(m-1)})=\sum_{i_1=1}^n\cdots\sum_{i_{(k+h)(m-1)}=1}^na_{i_1i_2}a_{i_2i_3}\cdots
a_{i_{(k+h)(m-1)-1}i_{(k+h)(m-1)}}a_{i_{(k+h)(m-1)}i_1}
\end{eqnarray*}
which results in
$(a_{ii})^{k(m-1)-s}(a_{ij})^{s}(a_{ji})^{s}(a_{jj})^{h(m-1)-s}$ has
and only has either the packaged elements of $i$ or the packaged
elements of $j$ if we count from the left most in the expression,
and is totally determined by the numbers of $a_{jj}$'s in the
packaged elements and the positions of the packaged elements in the
expression
\begin{eqnarray}\label{new-7}
a_{i_1i_2}a_{i_2i_3}\cdots
a_{i_{(k+h)(m-1)-1}i_{(k+h)(m-1)}}a_{i_{(k+h)(m-1)}i_1}.
\end{eqnarray}
So, the number of term
$(a_{ii})^{k(m-1)-s}(a_{ij})^{s}(a_{ji})^{s}(a_{jj})^{h(m-1)-s}$ in
the expansion
\begin{eqnarray*}
\mbox{Tr}(A^{(k+h)(m-1)})=\sum_{i_1=1}^n\cdots\sum_{i_{(k+h)(m-1)}=1}^na_{i_1i_2}a_{i_2i_3}\cdots
a_{i_{(k+h)(m-1)-1}i_{(k+h)(m-1)}}a_{i_{(k+h)(m-1)}i_1}
\end{eqnarray*}
is totally determined by the number of cases how the packaged
elements are arranged multiplying the number of cases of the
positions of the packaged elements in the expression \reff{new-7}.

In the following, we consider only the situation of packaged
elements of $i$. The other situation is similar. Note that there are
altogether $s$ packaged elements of $i$ in every expression
\reff{new-7} which results in
$(a_{ii})^{k(m-1)-s}(a_{ij})^{s}(a_{ji})^{s}(a_{jj})^{h(m-1)-s}$.

Firstly, note that there are $h(m-1)-s$ $a_{jj}$'s in
$(a_{ii})^{k(m-1)-s}(a_{ij})^{s}(a_{ji})^{s}(a_{jj})^{h(m-1)-s}$.
Then we have
\begin{eqnarray*}
\left(\begin{array}{c}h(m-1)-s+(s-1)\\s-1\end{array}\right)=\left(\begin{array}{c}h(m-1)-1\\s-1\end{array}\right)
\end{eqnarray*}
different cases of $s$ packaged elements of $i$ with every case
results in $(a_{ij})^{s}(a_{ji})^{s}(a_{jj})^{h(m-1)-s}$. Secondly,
for an arbitrary but fixed case of $s$ packaged elements of $i$ in
the first step, we get $k(m-1)$ ``mixed" elements consist of the $s$
packaged elements and the rest $k(m-1)-s$ $a_{ii}$'s. Consequently,
there are exactly
\begin{eqnarray*}
\left(\begin{array}{c}k(m-1)\\s\end{array}\right)
\end{eqnarray*}
cases of the expression \reff{new-7}, under which the expression
\reff{new-7} results in the term\\
$(a_{ii})^{k(m-1)-s}(a_{ij})^{s}(a_{ji})^{s}(a_{jj})^{h(m-1)-s}$.

Consequently, the number of term
$(a_{ii})^{k(m-1)-s}(a_{ij})^{s}(a_{ji})^{s}(a_{jj})^{k(m-1)-s}$ in
the expansion $\mbox{Tr}(A^{(k+h)(m-1)})$ is
\begin{eqnarray*}
\left(\begin{array}{c}k(m-1)\\s\end{array}\right)\left(\begin{array}{c}h(m-1)-1\\s-1\end{array}\right)
\end{eqnarray*}
in the situation of packaged elements of $i$.

By the symmetry of $i$ and $j$, we can prove similarly that, in the
situation of packaged elements of $j$, the number of term
$(a_{ii})^{k(m-1)-s}(a_{ij})^{s}(a_{ji})^{s}(a_{jj})^{k(m-1)-s}$ in
the expansion $\mbox{Tr}(A^{(k+h)(m-1)})$ is
\begin{eqnarray*}
\left(\begin{array}{c}h(m-1)\\s\end{array}\right)\left(\begin{array}{c}k(m-1)-1\\s-1\end{array}\right).
\end{eqnarray*}
Hence, we have that the number of term
$(a_{ii})^{k(m-1)-s}(a_{ij})^{s}(a_{ji})^{s}(a_{jj})^{k(m-1)-s}$ in
the expansion $\mbox{Tr}(A^{(k+h)(m-1)})$ is altogether
\begin{eqnarray*}
\left(\begin{array}{c}k(m-1)\\s\end{array}\right)\left(\begin{array}{c}h(m-1)-1\\s-1\end{array}\right)
+\left(\begin{array}{c}h(m-1)\\s\end{array}\right)\left(\begin{array}{c}k(m-1)-1\\s-1\end{array}\right).
\end{eqnarray*}
The proof is complete.  \ep

In the sequel, let $\hat g$ and $g$ be defined as in
\reff{cha-1-1}and \reff{cha-1-2}, respectively.
\begin{Lemma}\label{lem-2}
Let ${\cal T}\in \mathbb{T}(\mathbb{C}^n,m)$. For arbitrary $i<j$,
and $h,k\geq 1$, we have
\begin{eqnarray}\label{lem-2-1}
&&\frac{(\hat g_i)^h(\hat g_j)^k}{(h(m-1))!(k(m-1))!}\mbox{Tr}(A^{(h+k)(m-1)})\nonumber\\
&=&(\frac{h+k}{hk(m-1)})\sum_{s=1}^{\min\{h,k\}(m-1)}\sum_{(a_1,\ldots,a_h)\in\mathbb{D}^s,\;(b_1,\ldots,b_k)\in\mathbb{E}^s}s\displaystyle\prod_{p=1}^h\displaystyle\prod_{q=1}^kw_{ii\ldots
i\underbrace{j\ldots j}_{a_p}}w_{j\underbrace{i\ldots
i}_{b_q}j\ldots j},\;\;\;\;\;\;\;
\end{eqnarray}
with $\mathbb{D}^s:=\{(a_1,\ldots,a_h)\;|\;a_1+\cdots+a_h=s,\;0\leq
a_p\leq m-1\;\forall p=1,\ldots,h\}$ and
$\mathbb{E}^s:=\{(b_1,\ldots,b_k)\;|\;b_1+\cdots+b_k=s,\;0\leq
b_q\leq m-1\;\forall q=1,\ldots,k\}$.
\end{Lemma}

\noindent {\bf Proof.} Let $w:=\min\{h,k\}(m-1)$,
$\mathbb{D}^s:=\{(a_1,\ldots,a_h)\;|\;a_1+\cdots+a_h=s,\;0\leq
a_p\leq m-1\;\forall p=1,\ldots,h\}$ and
$\mathbb{E}^s:=\{(b_1,\ldots,b_k)\;|\;b_1+\cdots+b_k=s,\;0\leq
b_q\leq m-1\;\forall q=1,\ldots,k\}$ for all $s=1,\ldots,w$. By
\reff{cha-1-1} and \reff{prop-3-1}, we have
\begin{eqnarray*}
&&(\hat g_i)^h(\hat g_j)^k\mbox{Tr}(A^{(h+k)(m-1)})\\
&=&\left[\sum_{i_2\leq\cdots\leq i_m}w_{ii_2\ldots i_m}\frac{\partial}{\partial a_{ii_2}}\cdots \frac{\partial}{\partial a_{ii_m}}\right]^h\left[\sum_{j_2\leq\cdots\leq j_m}w_{jj_2\ldots j_m}\frac{\partial}{\partial a_{jj_2}}\cdots \frac{\partial}{\partial a_{jj_m}}\right]^k\mbox{Tr}(A^{(h+k)(m-1)})\\
&=&\left[\sum_{i_2\leq\cdots\leq i_m}w_{ii_2\ldots i_m}\frac{\partial}{\partial a_{ii_2}}\cdots \frac{\partial}{\partial a_{ii_m}}\right]^h\left[\sum_{j_2\leq\cdots\leq j_m}w_{jj_2\ldots j_m}\frac{\partial}{\partial a_{jj_2}}\cdots \frac{\partial}{\partial a_{jj_m}}\right]^k\\
&&\cdot\left(\sum_{i_1=1}^n\ldots\sum_{i_{(h+k)(m-1)}=1}^na_{i_1i_2}a_{i_2i_3}\cdots a_{i_{(h+k)(m-1)-1}i_{(h+k)(m-1)}}a_{i_{(h+k)(m-1)}i_1}\right)\\
&=&\left(\sum_{s=1}^{w}\sum_{(a_1,\ldots,a_h)\in\mathbb{D}^s,\;(b_1,\ldots,b_k)\in\mathbb{E}^s}\displaystyle\prod_{p=1}^h\displaystyle\prod_{q=1}^kw_{ii\ldots i\underbrace{j\ldots j}_{a_p}}w_{j\underbrace{i\ldots i}_{b_q}j\ldots j}\right.\\
&&\left.\cdot\left(\frac{\partial}{\partial
a_{ii}}\right)^{h(m-1)-s}
\left(\frac{\partial}{\partial a_{ij}}\right)^s\left(\frac{\partial}{\partial a_{ji}}\right)^s\left(\frac{\partial}{\partial a_{jj}}\right)^{k(m-1)-s}\right)\\
&&\cdot \left\{\sum_{s=1}^{w}\left[\left(\begin{array}{c}k(m-1)\\s\end{array}\right)\left(\begin{array}{c}h(m-1)-1\\s-1\end{array}\right)\right.\right.\\
&&+\left.\left.\left(\begin{array}{c}h(m-1)\\s\end{array}\right)\left(\begin{array}{c}k(m-1)-1\\s-1\end{array}\right)\right]\cdot(a_{ii})^{h(m-1)-s}(a_{ij})^{s}(a_{ji})^{s}(a_{jj})^{k(m-1)-s}\right\}\\
&=&\sum_{s=1}^{w}\sum_{(a_1,\ldots,a_h)\in\mathbb{D}^s,\;(b_1,\ldots,b_k)\in\mathbb{E}^s}\displaystyle\prod_{p=1}^h\displaystyle\prod_{q=1}^kw_{ii\ldots i\underbrace{j\ldots j}_{a_p}}w_{j\underbrace{i\ldots i}_{b_q}j\ldots j}\\
&&\;\;\cdot s\left((k(m-1))!(h(m-1)-1)!+(h(m-1))!(k(m-1)-1)!\right).
\end{eqnarray*}
Here, the third equality follows from Lemma \ref{lem-2-0}.
Consequently, \reff{lem-2-1} follows.  \ep

Especially, the following is a direct corollary of Lemma
\ref{lem-2}.

\begin{Corollary}\label{cor-2}
Let ${\cal T}\in \mathbb{T}(\mathbb{C}^n,m)$. For arbitrary $i<j$,
we have
\begin{eqnarray*}
\frac{\hat g_i\hat
g_j}{[(m-1)!]^2}\mbox{Tr}(A^{2(m-1)})=\sum_{s=1}^{m-1}\left(\frac{2s}{m-1}\right)w_{ii\ldots
i\underbrace{j\ldots j}_s}w_{ji\ldots i\underbrace{j\ldots
j}_{m-1-s}}.
\end{eqnarray*}
\end{Corollary}

Given an index set $\mathbb{L}:=\{k_1,\ldots,k_l\}$ with $k_s$
taking value in $\{1,\ldots,n\}$, denote by
$\mathbb{H}_i(\mathbb{L})$ the set of indices in $\mathbb{L}$ taking
value $i$.  We denote by $|\mathbb{E}|$ the cardinality of a set
$\mathbb{E}$.

\begin{Proposition}\label{prop-5}
Let ${\cal T}\in \mathbb{T}(\mathbb{C}^n,m)$. We have
\begin{eqnarray*}
\mbox{Tr}_2({\cal T})&=&(m-1)^{n-1}\left[\sum_{i=1}^n\frac{(\hat g_i)^2}{(2(m-1))!}+\sum_{i<j}\frac{\hat g_i\hat g_j}{[(m-1)!]^2}\right]\mbox{Tr}(A^{2(m-1)})\nonumber\\
&=&(m-1)^{n-1}\left[\sum_{i=1}^nt_{ii\ldots i}^2+\sum_{i<j}\sum_{s=1}^{m-1}\left(\frac{2s}{m-1}\right)\right.\nonumber\\
&&\;\;\;\left.\cdot \left(\sum_{|\mathbb{H}_i(\{i_2,\ldots,i_m\})|=m-1-s,\;|\mathbb{H}_j(\{i_2,\ldots,i_m\})|=s}t_{ii_2\ldots i_m}\right)\right.\nonumber\\
&&\;\;\;\;\cdot
\left.\left(\sum_{|\mathbb{H}_i(\{j_2,\ldots,j_m\})|=s,\;|\mathbb{H}_j(\{j_2,\ldots,j_m\})|=m-1-s}t_{jj_2\ldots
j_m}\right)\right].
\end{eqnarray*}
\end{Proposition}

\noindent {\bf Proof.} The result follows from Lemma \ref{lem-1},
Corollary \ref{cor-2}, \reff{cha-1} and \reff{cha-1-1} immediately.
\ep

When $n=2$, we can get the coefficients of the characteristic
polynomial $\psi(\lambda)$ explicitly in terms of the entries of the
underlying tensor by using Theorem \ref{thm-2}, Lemmas \ref{lem-1}
and \ref{lem-2}. It is an alternative to Sylvester's formula
\cite{sy}. In the following theorem, $w$, $\mathbb{D}^s$'s and
$\mathbb{E}^s$'s are defined as in those in Lemma \ref{lem-2}.

\begin{Theorem}\label{thm-3}
Let ${\cal T}\in \mathbb{T}(\mathbb{C}^2,m)$. We have
\begin{eqnarray*}
\psi(\lambda)&=&\lambda^{2(m-1)}+\sum_{k=1}^{2(m-1)}\lambda^{2(m-1)-k}\sum_{i=1}^k\frac{1}{i!}\sum_{d_j>0,\sum_{j=1}^id_j=k}\displaystyle\prod_{j=1}^i\frac{-\mbox{Tr}_{d_j}({\cal
T})}{d_j}
\end{eqnarray*}
with
\begin{eqnarray*}
&&\mbox{Tr}_d({\cal T})\nonumber\\
&=&(m-1)\left\{(t_{11\ldots 1}^d+t_{22\ldots 2}^d)+\sum_{h+k=d,h,k\geq 1}\sum_{s=1}^{w}\frac{s(h+k)}{hk(m-1)}\right.\nonumber\\
&&\cdot\left[\sum_{(a_1,\ldots,a_h)\in\mathbb{D}^s,\;(b_1,\ldots,b_k)\in\mathbb{E}^s}\displaystyle\prod_{p=1}^h\displaystyle\prod_{q=1}^k\left(\sum_{|\mathbb{H}_1(\{i_2,\ldots,i_m\})|=m-1-a_p,\;|\mathbb{H}_2(\{i_2,\ldots,i_m\})|=a_p}t_{1i_2\ldots i_m}\right)\right.\nonumber\\
&&\;\;\;\;\left.\left.\cdot
\left(\sum_{|\mathbb{H}_1(\{j_2,\ldots,j_m\})|=b_p,\;|\mathbb{H}_2(\{j_2,\ldots,j_m\})|=m-1-b_p}t_{2j_2\ldots
j_m}\right)\right]\right\}
\end{eqnarray*}
for $d\in\{1,\ldots,2(m-1)\}$.
\end{Theorem}

Note that as $\mbox{Edet}({\cal T})=(-1)^{n(m-1)^{n-1}}\psi(0)$ by
Theorem \ref{thm-2} and Corollary \ref{cor-001}, when $n=2$, we get
an explicit formula for $\mbox{Edet}({\cal T})$ as
\begin{eqnarray*}
\mbox{Edet}({\cal
T})=\sum_{i=1}^{2(m-1)}\frac{1}{i!}\sum_{d_j>0,\sum_{j=1}^id_j=2(m-1)}\displaystyle\prod_{j=1}^i\frac{-\mbox{Tr}_{d_j}({\cal
T})}{d_j}.
\end{eqnarray*}

When $m=3$, we have the following corollary.

\begin{Corollary}\label{cor-3}
Let ${\cal T}\in \mathbb{T}(\mathbb{C}^2,3)$. We have
\begin{eqnarray*}
\psi(\lambda)&=&\lambda^4-\lambda^3\mbox{Tr}_1({\cal T})+\frac{1}{2}\lambda^2\left(\left[\mbox{Tr}_1({\cal T})\right]^2-\mbox{Tr}_2({\cal T})\right)\\
&+&\frac{1}{12}\lambda\left(-2\left[\mbox{Tr}_1({\cal T})\right]^3+6\mbox{Tr}_1({\cal T})\mbox{Tr}_2({\cal T})-4\mbox{Tr}_3({\cal T})\right)\\
&+&\frac{1}{24}\left(\left[\mbox{Tr}_1({\cal
T})\right]^4-6\left[\mbox{Tr}_1({\cal T})\right]^2\mbox{Tr}_2({\cal
T})+8\mbox{Tr}_1({\cal T})\mbox{Tr}_3({\cal
T})+3\left[\mbox{Tr}_2({\cal T})\right]^2-6\mbox{Tr}_4({\cal
T})\right)
\end{eqnarray*}
and
\begin{eqnarray*}
\mbox{Edet}({\cal T})=\frac{1}{24}\left(\left[\mbox{Tr}_1({\cal
T})\right]^4-6\left[\mbox{Tr}_1({\cal T})\right]^2\mbox{Tr}_2({\cal
T})+8\mbox{Tr}_1({\cal T})\mbox{Tr}_3({\cal
T})+3\left[\mbox{Tr}_2({\cal T})\right]^2-6\mbox{Tr}_4({\cal
T})\right)
\end{eqnarray*}
with
\begin{eqnarray*}
\mbox{Tr}_1({\cal T})&=&2(t_{111}+t_{222}),\\
\mbox{Tr}_2({\cal T})&=&2(t_{111}^2+t_{222}^2)+2(t_{112}+t_{121})(t_{212}+t_{221})+4(t_{122}t_{211}),\\
\mbox{Tr}_3({\cal T})&=&2(t_{111}^3+t_{222}^3)\\
&&+\frac{3}{2}(t_{112}+t_{121})\left[(t_{212}+t_{221})t_{222}\right]+3t_{122}\left[(t_{212}+t_{221})^2+t_{211}t_{222}\right]\\
&&+\frac{3}{2}(t_{212}+t_{221})\left[(t_{112}+t_{121})t_{111}\right]+3t_{211}\left[(t_{112}+t_{121})^2+t_{122}t_{111}\right],\\
\mbox{Tr}_4({\cal T})&=&2(t_{111}^4+t_{222}^4)\\
&&+\frac{4}{3}(t_{112}+t_{121})\left[t_{222}^2(t_{212}+t_{221})\right]+\frac{8}{3}t_{122}\left[t_{222}(t_{212}+t_{221})^2+t_{222}^2t_{211}\right]\\
&&+\frac{4}{3}(t_{212}+t_{221})\left[t_{122}^2(t_{112}+t_{121})\right]+\frac{8}{3}t_{211}\left[t_{111}(t_{112}+t_{121})^2+t_{111}^2t_{122}\right]\\
&&+t_{111}(t_{121}+t_{121})t_{222}(t_{212}+t_{221})+3\left[t_{111}(t_{121}+t_{112})\right]\left[t_{211}(t_{221}+t_{212})\right]\\
&&\;\;+2\left[(t_{121}+t_{112})^2+t_{122}t_{111}\right]\left[t_{211}t_{222}+(t_{212}+t_{221})^2\right]+4t_{122}^2t_{211}^2.
\end{eqnarray*}
\end{Corollary}

\section{Inequalities of the E-Determinant}
 \setcounter{Assumption}{0}
\setcounter{Theorem}{0} \setcounter{Proposition}{0}
\setcounter{Corollary}{0} \setcounter{Lemma}{0}
\setcounter{Definition}{0} \setcounter{Remark}{0}
\setcounter{Algorithm}{0}  \setcounter{Example}{0} \hspace{4mm} In
this section, we generalize several inequalities of the determinants
for matrices to the E-determinants for tensors.

\subsection{Hadamard's Inequality}
\hspace{4mm} We generalize Hadamrd's inequality for matrices to
tensors in this subsection. We assume that $m$ is an even integer in
this subsection.

\begin{Lemma}\label{lem-6}
Let ${\cal T}\in \mathbb{S}(\mathbb{R}^n,m)$ and ${\cal U}$ be a sub-tensor of ${\cal T}$ associated to any nonempty subset of $\{1,\ldots,n\}$. If ${\cal T}$ is positive semidefinite, then $\mbox{Edet}({\cal U})\geq 0$.
\end{Lemma}
\noindent {\bf Proof.} It follows from the definitions of positive semidefiniteness and sub-tensors (Definition \ref{def-4}), and \cite[Proposition 7]{q05}. \ep

The following assumption is involved.
\begin{Assumption}\label{assump2}
Let ${\cal T}\in \mathbb{S}(\mathbb{R}^n,m)$. Suppose that
$\sigma({\cal T})$ consists of only real numbers.
\end{Assumption}

\begin{Lemma}\label{lem-3}
Let ${\cal T}\in \mathbb{S}(\mathbb{R}^n,m)$, and suppose that Assumptions \ref{assump} and \ref{assump2} hold. If ${\cal T}$ is positive semidefinite and $t_{i\ldots i}\leq 1$ for all $i\in\{1,\ldots,n\}$, then $1\geq\mbox{Edet}({\cal T})\geq 0$.
\end{Lemma}

\noindent {\bf Proof.} By Lemma \ref{lem-6}, $\mbox{Edet}({\cal T})\geq 0$. Furthermore,
\begin{eqnarray*}
n(m-1)^{n-1}&\geq& (m-1)^{n-1}\sum_{i=1}^nt_{ii\ldots i}\\
&=&\sum_{\lambda_i\in\sigma({\cal T})}m_i\lambda_i\\
&\geq& n(m-1)^{n-1}\left(\Pi_{\lambda_i\in\sigma({\cal T})}\lambda_i^{m_i}\right)^{\frac{1}{n(m-1)^{n-1}}}\\
&=&n(m-1)^{n-1}\left(\mbox{Edet}({\cal T})\right)^{\frac{1}{n(m-1)^{n-1}}},
\end{eqnarray*}
where $m_i$ is the algebraic multiplicity of eigenvalue $\lambda_i$. Here the first inequality follows from the assumption that $t_{i\ldots i}\leq 1$ for all $i\in\{1,\ldots,n\}$; the first equality from Corollary \ref{cor-1}; and, the last inequality from Assumptions \ref{assump} and \ref{assump2} (which, together with the positive semidefiniteness of ${\cal T}$, imply that ${\cal T}$ has only nonnegative eigenvalues) and the arithmetic geometry mean inequality. This, together with $\mbox{Edet}({\cal T})\geq 0$, implies $\mbox{Edet}({\cal T})\leq 1$. The proof is complete.  \ep

Note again that the tensor in \cite[Example 2]{q05} satisfies all the hypotheses in Lemma \ref{lem-3}.
\begin{Theorem}\label{thm-5}
Let ${\cal T}\in \mathbb{S}(\mathbb{R}^n,m)$,
$D:=\mbox{Diag}\{t_{1\ldots 1},\ldots,t_{n\ldots n}\}$ be the
diagonal matrix consisting of diagonal elements $t_{1\ldots
1},\ldots,t_{n\ldots n}$, and ${\cal K}:=D^{-\frac{1}{m}}\circ {\cal
T}\circ D^{-\frac{1}{m}}$, where $D$ is invertible. Suppose that
Assumptions \ref{assump} and \ref{assump2} hold for both ${\cal T}$
and ${\cal K}$. If ${\cal T}$ is positive semidefinite, then
$0\leq\mbox{Edet}({\cal T})\leq \left(\Pi_{i=1}^nt_{i\ldots
i}\right)^{(m-1)^{n-1}}$.
\end{Theorem}

\noindent {\bf Proof.} If $t_{i\ldots i}=0$ for some $i$,
then ${\cal T}\mathbf{e}_i^m:=\mathbf{e}_i^T({\cal T}\mathbf{e}_i^{m-1})=0$. Consequently, ${\cal T}$ has a zero eigenvalue, since it is positive semidefinite and symmetric (see proof of \cite[Theorem 5]{q05}). This, together with Corollary \ref{cor-1}, further implies $\mbox{Edet}({\cal T})=0$. Then, the results follow trivially.

Now, suppose that matrix $D$ is invertible. By the definition of positive semidefiniteness, it is clear that ${\cal K}$ is positive semidefinite as well. Since it also satisfies Assumptions \ref{assump} and \ref{assump2}, consequently,
$0\leq \mbox{Edet}({\cal K})\leq 1$ by Lemma \ref{lem-3}.
Furthermore, by Propositions \ref{prop-2-2} and \ref{prop-2-3}, it is easy to see that
\begin{eqnarray*}
\mbox{Edet}({\cal K})=\left(\mbox{Det}(D^{-\frac{1}{m}})\right)^{m(m-1)^{n-1}}\mbox{Edet}({\cal T}).
\end{eqnarray*}
Hence,
\begin{eqnarray*}
\mbox{Edet}({\cal T})=\left(\Pi_{i=1}^nt_{i\ldots i}\right)^{(m-1)^{n-1}}\mbox{Edet}({\cal K})\leq \left(\Pi_{i=1}^nt_{i\ldots i}\right)^{(m-1)^{n-1}}.
\end{eqnarray*}

The proof is complete. \ep

\begin{Remark}\label{rmk-3}
{\em When $m=2$, Assumptions \ref{assump} and \ref{assump2} are satisfied by matrices, so Theorem \ref{thm-5} reduces to the well-known Hadamard's inequality \cite[Theorem 2.5.4]{hj2}}.
\end{Remark}

\subsection{Ger\v{s}gorin's Inequality}
\hspace{4mm} We generalize Ger\v{s}gorin's inequality for matrices \cite[Problem 6.1.3]{hj1} to tensors in this subsection.

\begin{Lemma}\label{lem-4}
Let ${\cal T}\in\mathbb{T}(\mathbb{C}^n,m)$ and $\rho({\cal T}):=\max_{\lambda\in\sigma({\cal T})}|\lambda|$ be its spectral radius. Then,
\begin{eqnarray*}
\rho({\cal T})\leq \max_{1\leq i\leq n}\left(\sum_{i_2,\ldots,i_m=1}^n|t_{ii_2\ldots i_m}|\right).
\end{eqnarray*}
\end{Lemma}

\noindent {\bf Proof.} The result follows from \cite[Theorem 6]{q05} immediately. \ep

\begin{Proposition}\label{prop-ge}
Let ${\cal T}\in\mathbb{T}(\mathbb{C}^n,m)$. Then,
\begin{eqnarray}\label{ga-2}
|\mbox{Edet}({\cal T})|\leq \displaystyle\prod_{1\leq i\leq n}\left(\sum_{i_2,\ldots,i_m=1}^n|t_{ii_2\ldots i_m}|\right)^{(m-1)^{n-1}}.
\end{eqnarray}
\end{Proposition}

\noindent {\bf Proof.} If $\sum_{i_2,\ldots,i_m=1}^n|t_{ii_2\ldots i_m}|=0$ for some $i\in\{1,\ldots,n\}$, then $\mbox{Edet}({\cal T})=0$ by Proposition \ref{prop-000} (i). Consequently, \reff{ga-2} follows trivially.

Now, suppose that $\sum_{i_2,\ldots,i_m=1}^n|t_{ii_2\ldots i_m}|\neq 0$ for all $i\in\{1,\ldots,n\}$. Let tensor ${\cal U}\in\mathbb{T}(\mathbb{C}^n,m)$ be defined as
\begin{eqnarray}\label{ga-3}
u_{ii_2\ldots i_m}:=\frac{t_{ii_2\ldots i_m}}{\sum_{i_2,\ldots,i_m=1}^n|t_{ii_2\ldots i_m}|},\;\forall i,i_2,\ldots,i_m=1,\ldots,n.
\end{eqnarray}
Then, by Lemma \ref{lem-4}, we have that $\rho({\cal U})\leq 1$. This, together with Corollary \ref{cor-1}, further implies that
\begin{eqnarray*}
|\mbox{Edet}({\cal U})|\leq 1.
\end{eqnarray*}
Moreover, by Proposition \ref{prop-000}(ii) and \reff{ga-3}, it is clear that
\begin{eqnarray*}\label{ga-4}
|\mbox{Edet}({\cal U})|=\frac{|\mbox{Edet}({\cal T})|}{\displaystyle\prod_{1\leq i\leq n}\left(\sum_{i_2,\ldots,i_m=1}^n|t_{ii_2\ldots i_m}|\right)^{(m-1)^{n-1}}}.
\end{eqnarray*}
Consequently, \reff{ga-2} follows. The proof is complete. \ep

\subsection{Minikowski's Inequality}
\hspace{4mm} We give a partial generalization of Minikowski's
inequality for matrices \cite[Theorem 7.8.8]{hj1} to tensors in this
subsection. We present the following lemma first.

\begin{Lemma}\label{lem-5}
Let ${\cal T}\in \mathbb{T}(\mathbb{C}^n,m)$. Then,
\begin{eqnarray*}
\lambda\in\sigma({\cal T})\quad\Longleftrightarrow\quad 1+\lambda\in\sigma({\cal E}+{\cal T}).
\end{eqnarray*}
\end{Lemma}

\noindent {\bf Proof.} It follows from Definition \ref{def-1} directly. \ep

\begin{Proposition}\label{prop-min}
Let ${\cal T}\in \mathbb{S}(\mathbb{R}^n,m)$ and $m$ be even, and suppose that Assumptions \ref{assump} and \ref{assump2} hold. If ${\cal T}$ is positive semidefinite, then
\begin{eqnarray*}
\left[\mbox{Edet}({\cal E}+{\cal T})\right]^{\frac{1}{n(m-1)^{n-1}}}\geq 1+\left[\mbox{Edet}({\cal T})\right]^{\frac{1}{n(m-1)^{n-1}}}.
\end{eqnarray*}
\end{Proposition}

\noindent {\bf Proof.} Suppose that $0\leq \lambda_1\leq\ldots\leq\lambda_{n(m-1)^{n-1}}$ are the eigenvalues of tensor ${\cal T}$. Then,
the arithmetic geometric mean inequality implies that
\begin{eqnarray*}
\Pi_{i=1}^{n(m-1)^{n-1}}(1+\lambda_i)\geq\left(1+\sqrt[n(m-1)^{n-1}]{\Pi_{i=1}^{n(m-1)^{n-1}}\lambda_i}\right)^{n(m-1)^{n-1}}.
\end{eqnarray*}
This, together with Corollary \ref{cor-1} and Lemma \ref{lem-5}, implies the result.  \ep


\section{Final Remarks}
 \setcounter{Assumption}{0}
\setcounter{Theorem}{0} \setcounter{Proposition}{0}
\setcounter{Corollary}{0} \setcounter{Lemma}{0}
\setcounter{Definition}{0} \setcounter{Remark}{0}
\setcounter{Algorithm}{0}  \setcounter{Example}{0} \hspace{4mm} In
this paper, we introduced the E-determinant of a tensor and
investigated its various properties. 
The simple application in Section 6 demonstrates that the
E-determinant theory is applicable and worth further exploring.

Certainly, it is worth exploring if Conjecture \ref{conj-00} is true
or not and what is the right expression of the E-determinant of the
composition if the conjecture is not true.   There are also many
other issues for further research in the theory of the
E-determinant, which is surely one of the foundations of the
eigenvalue theory of tensors. For example,
\begin{itemize}
\item How to remove Assumption \ref{assump} in Proposition \ref{prop-4-1}, Assumptions \ref{assump} and \ref{assump2} in Proposition \ref{prop-min}, Lemma \ref{lem-3} and Theorem \ref{thm-5}, or to derive which class of tensors satisfies these assumptions.
\item More explicit formulae for the characteristic polynomials of general tensors, like those in Theorem \ref{thm-3} and Corollary \ref{cor-3}.
\item More properties about the E-determinants. For example, given a matrix $A$, Laplace's formula \cite[Page 7]{hj1} reads:
\begin{eqnarray}\label{rm-1}
\mbox{Det}(A)=\sum_{j=1}^n(-1)^{i+j}a_{ij}M_{ij}
\end{eqnarray}
with minor $M_{ij}$ being defined to be the determinant of the
$(n-1)\times(n-1)$ matrix that results from $A$ by removing the
$i$-th row and the $j$-th column. If  a generalization of
\reff{rm-1} can be derived for the E-determinant of a tensor, many
other useful inequalities, for example, Oppenheim's inequality, can
be proved for the E-determinants.
\end{itemize}



\end{document}